\documentclass[12pt]{amsart}

\usepackage[a4paper,top=2.5cm,bottom=2.60cm,left=2.8cm,right=2.8cm]{geometry} 

\usepackage{amsmath}
\usepackage{amsfonts}
\usepackage{amsbsy}
\usepackage{amssymb}
\usepackage[normalem]{ulem}
\usepackage{times}
\usepackage{mathrsfs}
\usepackage{exscale}
\usepackage{graphicx}

\allowdisplaybreaks

\usepackage[colorlinks,citecolor=blue,linkcolor=blue,
            bookmarksopen,
            bookmarksnumbered
           ]{hyperref}
           
\usepackage{tikz}
\usetikzlibrary{patterns}

\usepackage{times}

\usepackage{color}
\definecolor{gree}{rgb}   {0.,   0.66,   0.25 }
\definecolor{mg}{rgb}   {0.9,  0.,    0.9}
\definecolor{marin}{rgb}   {0.,   0.65,   0.5}
\definecolor{orange}{rgb}   {1.,   0.6,   0.}
\definecolor{brown}{rgb}   {0.4,   0.,   0.8}
\definecolor{grey}{rgb}   {0.6,   0.6,   0.6}
\definecolor{greymg}{rgb}   {1,   0.,   0.8}
\definecolor{greygr}{rgb}   {0.5,   0.8,   0.5}

\newtheorem{theorem}{Theorem}[section]
\newtheorem{lemma}[theorem]{Lemma}
\newtheorem{proposition}[theorem]{Proposition}

\theoremstyle{definition}

\theoremstyle{remark}
\newtheorem{remark}[theorem]{Remark}
\newtheorem{example}[theorem]{Example}

\numberwithin{equation}{section}

\newcommand{\ee}{\hskip0.15ex}
\newcommand{\me}{\hskip-0.15ex}
\newcommand{\dd}[1]{_{\raise-0.6ex\hbox{$\scriptstyle #1$}}}

\newcommand {\Norm}[2]{ \mathchoice
    {|\ee #1\ee|\dd{#2}\,}
    {| #1 |_{#2}}
    {| #1 |_{#2}}
    {| #1 |_{#2}} }
\newcommand {\DNorm}[2]{ \mathchoice
    {\|\ee #1\ee\|\dd{#2}\,}
    {\| #1 \|_{#2}}
    {\| #1 \|_{#2}}
    {\| #1 \|_{#2}} }
\newcommand {\Normc}[3]{ \mathchoice
    {|\ee #1\ee|\dd{#2}^{#3}}
    {| #1 |_{#2}^{#3}}
    {| #1 |_{#2}^{#3}}
    {| #1 |_{#2}^{#3}} }
\newcommand {\DNormc}[3]{ \mathchoice
    {\|\ee #1\ee\|\dd{#2}^{#3}}
    {\| #1 \|_{#2}^{#3}}
    {\| #1 \|_{#2}^{#3}}
    {\| #1 \|_{#2}^{#3}} }

\newcommand\C{{\mathbb C}}
\newcommand\D{{\mathbb D}}
\newcommand\R{{\mathbb R}}
\newcommand\N{{\mathbb N}}

\newcommand\T{{\mathbb T}}
\newcommand\Z{{\mathbb Z}}

\newcommand{\sC}{{\mathscr C}}

\newcommand{\sj}{^{(j)}}

\newcommand\cI{{\mathcal{I}}}
\newcommand\cH{{\mathcal{H}}}

\newcommand\cK{{\mathcal{K}}}

\newcommand\cO{{\mathcal{O}}}

\newcommand\cS{{\mathcal{S}}}

\newcommand{\re}{{\mathrm e}}

\newcommand {\gU}{{\mathfrak U}}

\newcommand {\gu}{{\mathfrak u}}

\newcommand{\Ex}{{\sf Ex}}
\newcommand{\BL}{{\sf BL\; }}

\usepackage{soul}

\def\ED{\end{document}}

\begin{document}

\title[End-point non-regularity of the Dirichlet problem on Lipschitz domains]{End-point non-regularity of the Dirichlet problem on Lipschitz domains -- An elementary proof}

\author{Martin Costabel}

\email{martin.costabel@univ-rennes.fr}
\urladdr{https://costabel.perso.math.cnrs.fr}

\author{Monique Dauge}

\email{monique.dauge@univ-rennes.fr}
\urladdr{https://dauge.pages.math.cnrs.fr}

\address{IRMAR UMR 6625 du CNRS, Universit\'{e} de Rennes}
\address{Campus de Beaulieu,
35042 Rennes Cedex, France \bigskip}

\date{\bf\today}

\keywords{Dirichlet problem, Lipschitz domain, endpoint regularity, 
Sobolev norm, boundary layer}

\subjclass{35B65, 35J25, 46E35}

\begin{abstract}
We construct a plane Lipschitz domain for which the range of the Laplace operator from $H^{3/2}\cap H^1_0$ to $H^{-1/2}$ is not closed. 
We show that rapid oscillations of the boundary create a boundary layer that leads to unbounded $H^{3/2}$ norm of the solution of the Dirichlet problem. 
Our proof is elementary in the sense that it is directly based on homogeneity arguments and does not use deep results of harmonic analysis such as estimates of area integrals or harmonic measures of Dahlberg or Jerison and Kenig.
\end{abstract}

\maketitle

{\footnotesize
\parskip 1pt
\tableofcontents
}

\section{Introduction}

\subsection{Motivation and main results}

Let $\cO$ be a bounded domain in $\R^n$. The Dirichlet problem for the Laplacian $\Delta$ is realized as the operator 
\[
   \Delta:\quad H^1_0(\cO) \to H^{-1}(\cO)
\]
which is an isomorphism. Here $H^1_0(\cO)$ is the closure in the Sobolev space $H^1(\cO)$ of smooth functions with compact support in $\cO$, and $H^{-1}(\cO)$ is its dual.

If the boundary of $\cO$ is smooth, then any improvement of regularity of right hand sides implies a regularity shift for solutions: In the scale of Hilbert Sobolev spaces, this means that for any positive $s$
\[
   \Delta:\quad H^{s+1}\cap H^1_0(\cO) \to H^{s-1}(\cO)\quad
   \mbox{is an isomorphism}.
\]
This ``regularity shift'' property is shared by the much wider class of elliptic boundary value problems, see for instance Agmon-Douglis-Nirenberg \cite{AgmonDouglisNirenberg59}, Lions-Magenes \cite{LionsMagenes68}. 

If the boundary of $\cO$ is {\em not} smooth, the regularity shift does not persist as is, but has to be modified according to the assumptions on the boundary $\partial\cO$. A great variety of assumptions on domains $\cO$ have been considered in the literature. Let us mention:
\smallskip

{\em (a) Domains with conical singularities}, see Kondrat'ev \cite{Kondratev67} and Maz'ya-Plamenevskii \cite{MazyaPlamenevskii84a}, which includes polygons in dimension $2$, see Grisvard \cite{Grisvard85}, and also \cite{DBook}: For any $s>0$, except in a discrete set $\cS=\{s_1<s_2<\cdots\}$ with $s_1>0$,
\[
   \Delta:\quad H^{s+1}\cap H^1_0(\cO) \to H^{s-1}(\cO)\quad
   \mbox{is Fredholm}.
\]
This operator is still an isomorphism if $s\in[0,s_1)$, but is not surjective when $s\ge s_1$. On a polygon $\cO$, we have a simple formula for the forbidden set $\cS$: with $\omega_j$ the openings of the angles of $\cO$
\[
   \cS = \bigcup_j \big\{\tfrac{k\pi}{\omega_j},\quad k\in\N_1\big\}.
\]
 When a domain with conical singularities is moreover Lipschitz, the characterization of $s_1$ given in the above references yields that in any dimension, $s_1>\frac{1}{2}$.

{\em (b) Corner domains}, see Maz'ya et al \cite{MazyaPlamenevskii77,MazyaRossmann10}, and \cite{DBook}. In dimension $n\ge3$, such domains may have edges besides corners. There exist characterizations of Fredholm and semi-Fredholm properties for $\Delta$. A remarkable result holds if $\cO$ is a polyhedron (i.e. its boundary is a finite union of flat faces---segments in dimension 2---
\[
   \Delta:\quad H^{2}\cap H^1_0(\cO) \to L^2(\cO)\quad
   \mbox{has a closed range}.
\]
In other words, this operator is semi-Fredholm. In general it is not surjective (as soon as $\cO$ has non-convex corners ($n=2$) or edges ($n\ge3$).
\smallskip

{\em (c) Lipschitz domains}. 
Note that families of domains {\em (a)} or {\em (b)} may include non-Lipschitz domains: In the simplest situation of polygons, this happens if some angles have the opening $\omega=2\pi$, which corresponds to cracks, and then $s_1=\frac{1}{2}$.
If the polygon is Lipschitz, then $s_{1}-\frac12$ is positive, but can be arbitrarily small. 
\smallskip

This rapid description of the landscape shows why $s=\frac{1}{2}$ is an {\em end-point for Lipschitz regularity} and makes it natural to ask what happens for this value of $s$.
\smallskip

The aim of this paper is to give an elementary and self-contained 
proof of the following result.
\begin{theorem}
 \label{T:NotH32}
There exists a bounded Lipschitz domain $\cO\subset\R^{2}$ such that
$$
 \Delta : \quad H^{\frac32}(\cO)\cap H^{1}_{0}(\cO) \,\to\, H^{-\frac12}(\cO)
 \quad \mbox{has a non-closed range}.
$$
In particular it is not surjective.
\end{theorem}

Let us clarify that $H^{-\frac12}(\Omega)$ is the dual space of the subspace, denoted by $\tilde H^{\frac12}(\Omega)$ in the tradition of Triebel or $H^{\frac12}_{00}(\Omega)$ in the tradition of Lions-Magenes, of the Sobolev space $H^{\frac12}(\Omega)$ of functions whose extension by zero outside $\Omega$ belongs to $H^{\frac12}(\R^{2})$.
\smallskip

This theorem is due to Jerison and Kenig \cite[Theorem 0.4]{JerisonKenig1995}. Their proof uses deep tools from harmonic analysis on Lipschitz domains, such as 
estimates for harmonic measures and layer potentials, and in particular Dahlberg's area integral estimates \cite{Dahlberg1980}. 
These estimates imply that any harmonic function in $H^{\frac32}(\Omega)$ on a bounded Lipschitz domain $\Omega$ has a trace on the boundary in $H^{1}(\partial\Omega)$. 
On the other hand, a simple and explicit construction due to Guy David 
\cite[p.~176]{JerisonKenig1995}
gives an example of a function $v\in H^{\frac32}(\Omega)$ whose boundary trace $g$ does not belong to $H^{1}(\partial\Omega)$.
Then the unique solution $w\in H^{1}_{0}(\Omega)$ of 
$\Delta w=-\Delta v\in H^{-\frac12}(\Omega)$ cannot belong to $H^{\frac32}(\Omega)$, because the harmonic function $u=v+w$ with boundary trace $g$ cannot belong to $H^{\frac32}(\Omega)$.
\smallskip

Our motivation to contribute a new, more direct, proof of Theorem \ref{T:NotH32}, comes from the posting of the three prepublications by Amrouche and Moussaoui \cite{AM2023,AM2025,AM2026}, which try to question the validity of Dahlberg's result by ``proving'' an anti-Theorem \ref{T:NotH32} (according to \cite[Theorem 7.1]{AM2025} the Laplacian would be an isomorphism from $H^{\frac32}\cap H^{1}_{0}(\cO)$ onto $H^{-\frac12}(\cO)$). In their recent  preprint \cite{JK2026} Jerison and Kenig argue that the proof of the anti-Theorem suffers from a gap due to non-controlled interpolation norms in varying domains. In the paper \cite[\S 5.3.2 ]{CoDa2026}, we have already described a neighboring situation (polygons with approaching corners) in which interpolation norms do have a non-uniform behavior with respect to domain variation. 
In addition, \cite{JK2026} thoroughly revisits Dahlberg's result and different proofs of it.

Our contribution to this controversy is to show how an elementary homogeneity argument adapted to Guy David's construction allows to bypass the estimates of norms of boundary traces and to estimate the $H^{\frac32}$ norm of $w$ directly.
This method of proof does not use the above-mentioned deep results from harmonic analysis.

The standard closed graph theorem (or Banach open mapping theorem), combined with the fact that 
$$
 \Delta : \quad H^{\frac32}(\cO)\cap H^{1}_{0}(\cO) \,\to\, H^{-\frac12}(\cO)
 \quad \mbox{ is continuous and injective},
$$
implies that the non-closedness of its range is equivalent to the non-existence of an a-priori estimate
\begin{equation}
\label{E:apriori}
 \DNorm{u}{H^{\frac32}(\cO)} \le C\, \DNorm{\Delta u}{H^{-\frac12}(\cO)} 
 \quad\mbox{ for all }u\in H^{\frac32}(\cO)\cap H^{1}_{0}(\cO) \,,
\end{equation}
or again to the existence of a sequence $(u_{j})_{j\in\N}$ 
in $H^{\frac32}(\cO)\cap H^{1}_{0}(\cO)$
such that $\Delta u_{j}$ remains bounded in $H^{-\frac12}(\cO)$, whereas
$$
 \DNorm{u_{j}}{H^{\frac32}(\cO)} \to\infty
 \quad\mbox{ as }\; j\to\infty \,.
$$

In any case, the example with infinite $H^{\frac32}$ norm is constructed from gluing together a sequence of domains for which one has finite, but not uniformly bounded, $H^{\frac32}$ norms. Thus, the main point is to find examples for the following situation.

\begin{theorem}
 \label{T:blowup}
There exist a plane domain $\cO_0$ and a sequence of domains $(\cO_{N})_{N\ge1}$ such that
\begin{equation}
\label{eq:ON}
\begin{cases}
   \cO_N\subset\cO_0 & \forall N\ge1\\[0.5ex]
   \partial\cO_N \underset{N\to\infty}{\to} \partial\cO_0\quad&\mbox{in Hausdorff distance}\\
   (\partial\cO_N)_{N\ge1} \quad&\mbox{are uniformly Lipschitz}
\end{cases}
\end{equation} 
and satisfying the two equivalent conditions (i) and (ii):
\smallskip

{\rm (i)}
There exists a bounded sequence $(v_{N})_{N\ge1}$, in the Sobolev space 
$H^{\frac32}(\cO_{0})$ such that 
 with  the trace $g_{N}$ of $v_{N}$ on $\Gamma_{N}=\partial\cO_{N}$ and 
 the solution $u_{N}\in H^{1}(\cO_{N})$ of the Dirichlet problem 
\begin{equation}
\label{E:DirN}
  \Delta u_{N} =0 \quad\mbox{ in } \;\cO_N\,,\qquad
  u_{N} = g_{N} \quad\mbox{ on }\;\Gamma_{N}\,,
\end{equation}
we have $u_{N}\in H^{\frac32}(\cO_{N})$, but
\begin{equation}
\label{E:uNunbd}
 \DNorm{u_{N}}{H^{\frac32}(\cO_{N})} \to\infty
 \quad\mbox{ as }\; N\to\infty.
\end{equation}
\smallskip

{\rm (ii)}
There exists a bounded sequence $(f_{N})_{N}$ in 
$H^{-\frac12}(\cO_{0})$ such that for the solution $w_{N}$ of the  Dirichlet problem
\begin{equation}
\label{E:DirfN}
 \Delta w_{N} = f_{N}\quad\mbox{ in } \;\cO_{N} , \quad
  w_{N} = 0 \quad\mbox{ on }\;\Gamma_{N}\,,
\end{equation}
 there holds $w_{N} \in H^{\frac32}(\cO_{N})$
 and $\DNorm{w_{N}}{H^{\frac32}(\cO_{N})} \to\infty$ \;as\; $N\to\infty$.
\end{theorem}

\begin{remark}
In this work, we use the characterization of Sobolev spaces of positive order by simple or double integrals (see sec. \ref{SS:scal}). This makes elementary the proof of {\em (i)-(ii)}:
\smallskip

{\em (i)} If the sequence $\DNorm{v_{N}}{H^{\frac32}(\cO_{0})}$ is bounded, then, by restriction, the sequence $\DNorm{v_{N}}{H^{\frac32}(\cO_{N})}$ is bounded (by the same bound).
\smallskip

{\em (ii)} If the sequence $\DNorm{f_{N}}{H^{-\frac12}(\cO_{0})}$ is bounded, the functionals $f_N$ define ---by the adjoint of the extension by zero from $\cO_N$ to $\cO_0$--- a bounded sequence $\DNorm{f_{N}}{H^{-\frac12}(\cO_{N})}$ (with the same bound).
\smallskip

{\em (iii)} The equivalence between the two formulations \eqref{E:DirN} and \eqref{E:DirfN} of the inhomogeneous Dirichlet problem is then given by the relations
$$
w_{N}=u_{N}-v_{N}\quad\mbox{ and } \quad f_{N}=-\Delta v_{N}\,,
$$
combined with continuity of $\Delta:H^{\frac32}(\cO_0)\to H^{-\frac12}(\cO_0)$.
\smallskip

{\em (iv)} Whereas we use the traces $g_{N}$ in our construction, the arguments do not rely on any estimates of their norms. However, it can be seen that the $g_N$ that we will construct are such that their $H^1(\Gamma_N)$-norms blow up as $N\to\infty$.
\end{remark}

\subsection{Explicit families of Lipschitz domains $\cO_N$}
We will prove that the conclusions of Theorem \ref{T:blowup} are satisfied for three distinct, but closely related, families of sequences $(\cO_{N})_{N\ge1}$, which we may call (a) {\em circular saws}, (b) {\em periodic saws}, and (c) {\em straight saws}. Their common feature is a generating motif (the tooth) that is scaled and repeated. Each one of these families allows to simplify some parts of the arguments, so that their simultaneous consideration serves the goal of maximal simplicity and transparency.
\smallskip

(a) \emph{\textbf {Circular saws}}: the limiting domain $\cO_0$ is the unit disk $\D$, the generating motif is $\cO_1=\Omega_1$ with boundary given in polar coordinates $(r,\theta)$ by $\theta\mapsto r=\rho(\theta)$ where $\rho$ is an arbitrarily chosen $2\pi$-periodic Lipschitz function, which we assume to be piecewise $\sC^2$ and to satisfy $R_0\le\rho(\theta)\le1$ with some $R_0>0$. The boundary of $\cO_N=\Omega_N$ is then given in polar coordinates by
$r=\rho_{N}(\theta)$ with
\[
   \rho_{N}(\theta)=\rho(N\theta)^{\frac1N}.
\]
\begin{proposition}
\label{P:a}
The sequence $(\cO_N)=(\Omega_N)$ converging to $\cO_0=\D$ is such that \eqref{eq:ON} holds. If $\rho$ is not constant, this sequence satisfies the equivalent conditions (i) and (ii) in Theorem \ref{T:blowup}.
\end{proposition}
\smallskip

(b)-(c) \emph{\textbf {Periodic/straight saws}}. In the $(x,y)$ plane we set (with the torus $T=\R/(2\pi\Z)$)
\[
   Q_{0,\rm per} = \T \times (0,\infty)\quad\mbox{and}\quad
   Q_{0} = (0,2\pi) \times (0,\infty).
\]
Let
\begin{equation}
\label{eq:h}
   \mbox{$h:x\mapsto y=h(x)$ \quad $2\pi$-periodic, Lipschitz and piecewise $\sC^2$, $\quad h\ge0$. }
\end{equation} 
With
\[
   h_N(x) = \frac{1}{N} h(Nx)
\]
set
\[
   Q_{N,\rm per} = \{(x,y)\in Q_{0,\rm per},\;\; y>h_N(x)\}\quad\mbox{and}\quad
   Q_{N} = \{(x,y)\in Q_{0},\;\; y>h_N(x)\}.
\]

We define associated families of bounded Lipschitz domains by setting, for some chosen $Y_0>\max h$:
\begin{equation}
\label{eq:wtQN}
   \widetilde Q_{N,\rm per} = \{(x,y)\in Q_{N,\rm per},\;\; y<Y_0\}\quad\mbox{and}\quad
   \widetilde Q_{N} = \{(x,y)\in Q_{N},\;\; y<Y_0\}.
\end{equation}
Let us clarify that the boundary of $\widetilde Q_{N,\rm per}$ has two connected components defined by the equations $y=h_N(x)$ and $y=Y_0$, respectively, while the boundary of $\widetilde Q_{N}$ has a single connected component (note that $\widetilde Q_{0}$ is the rectangle $(0,2\pi)\times(0,Y_0)$). As in case (a), there holds

\begin{proposition}
\label{P:bc}
The sequence $(\cO_N)=(\widetilde Q_{N,\rm per})$ converging to $\cO_0=\widetilde Q_{0,\rm per}$ is such that \eqref{eq:ON} holds. If the function $h$ is not constant, this sequence satisfies the equivalent conditions (i) and (ii) in Theorem \ref{T:blowup}. The same holds for the straight versions $\widetilde Q_{N}$ of the ``saws''.
\end{proposition}

\begin{remark}
The proofs of Propositions \ref{P:a} and \ref{P:bc} are closely linked with each other because $Q_{N,\rm per}$ and $\Omega_N$ are related by the conformal mapping
\begin{equation}
\label{E:xi}
 \xi: Q_{N,\rm per}\ni(x,y)\mapsto r\re^{i\theta}\in\Omega_N 
  \quad\mbox{with } \; r=\re^{-y}, \; \theta=x\,,
  \quad\mbox{i.\ e. }\; \xi(x,y)=\re^{i(x+iy)}\,,
\end{equation}
if the generating motifs $h$ and $\rho$ satisfy $\rho=\re^{-h}$, i.e. $h=-\log\rho$. 
\end{remark}

\begin{figure}[h]
\centering
\includegraphics[width=0.48\textwidth]{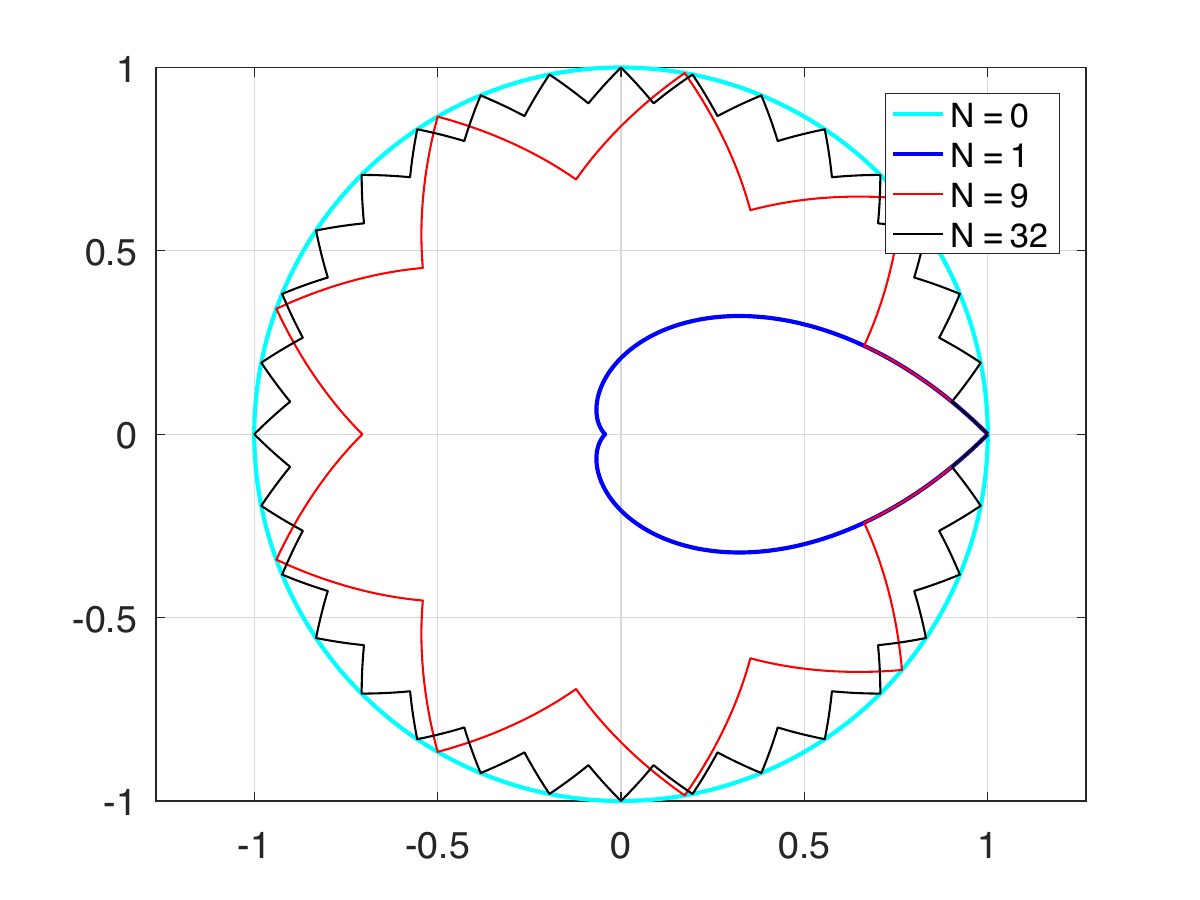}\quad
\includegraphics[width=0.48\textwidth]{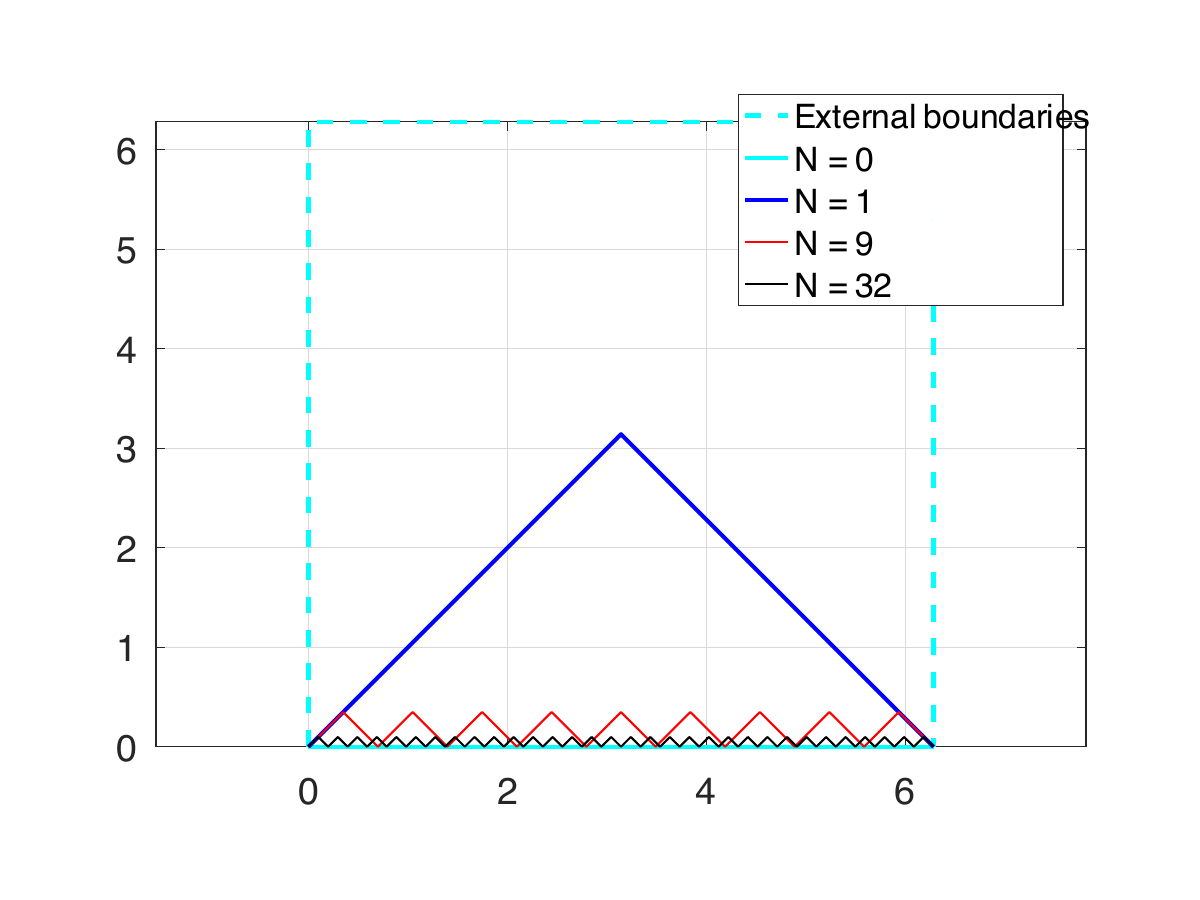}
\caption{\small Domains $\Omega_{N}$ (left) and $Q_N$ (right) for $N=1,9,32$ when the motif $h$ is the hat function \eqref{eq:hat} and $\rho=\re^{-h}$. The dashed vertical lines are boundaries of $Q_N$, but not of $Q_{N,\rm per}$. The dashed horizontal line is a boundary of $\widetilde Q_N$ and $\widetilde Q_{N,\rm per}$}
\label{F:hat}
\end{figure}

\begin{remark}
(i) The intuitive representation for the ``tooth'' defined by $h$ is a hat function:
\begin{equation}
\label{eq:hat}
   h(x) = x \;\mbox{ if }\; x\in(0,\pi) \quad\mbox{and}\quad 
   h(x) = 2\pi- x \;\mbox{ if }\; x\in(\pi,2\pi),
\end{equation}
see Figure \ref{F:hat}.
\smallskip

\noindent
(ii) Using $h$ defined by such a hat function is inspired by the example in \cite{JerisonKenig1995}, ascribed to Guy David. A sequence of domains congruent to $(\widetilde Q_{N})$ is used by \cite{AM2023}-\cite{AM2026} in the construction of their example meant to challenge Dahlberg's result.
\smallskip

\noindent
(iii) Our Proposition \ref{P:bc} includes the situation where $h$ is a smooth periodic function, in which case each Dirichlet problem $\Delta:H^{s+1}\cap H^1_0(\widetilde Q_{N,\rm per}) \to H^{s-1}(\widetilde Q_{N,\rm per})$ defines an isomorphism for any $s\ge0$. We stress that the accumulation of self-similar oscillations (regardless of their own regularity) is the cause for the $H^{\frac32}$ blow-up as $N\to\infty$, and not the presence of corners.
\end{remark}

\subsection{Sketch of proofs and plan of the paper}
To prove Proposition \ref{P:a}, we construct explicit sequences $(v_N)$ in $H^{\frac32}(\Omega_0)$ such that condition (i) in Theorem \ref{T:blowup} holds. The construction is performed in sec.\ \ref{S:disk} and estimates for the solution $u_N$ of \eqref{E:DirN} are proved in sec.\  \ref{S:est}. The proof of Proposition \ref{P:bc} stems from estimates obtained in $Q_{N,\rm per}$ via the mapping $\xi$ and from the use of cut-off functions, sec. \ref{S:Rect}. For a clue on the structure of $v_N$ and $u_N$, it is convenient to consider the corresponding cartesian variable functions $V_N$ and $U_N$ defined in $Q_{N,\rm per}$ (when $h=-\log\rho$) connected to $v_N$ and $u_N$ by the relations:
\[
   v_N(z) = V_N(\theta,-|\log z|) \quad\mbox{and}\quad u_N(z) = U_N(\theta,-|\log z|),
\]
with $z = |z|\re^{i\theta}$. The function $V_N$ depends only on the second variable and is defined as
\begin{equation}
\label{eq:deltaN}
   V_N(x,y) = A_{\delta_N}(y)\quad\mbox{with}\quad \delta_N = 
   \frac{|\log(\min\rho)|}{N} = \frac{\max h}{N}
\end{equation}
where the function $A_{\delta}(y)$ is linear, given by $\log\me|\log\delta|\;y$ when $0\le y\le \delta$, and is defined as a primitive of the function $t\mapsto \log|\log t|$ for larger values of $y$. The bounded harmonic solution $U_N$ in $Q_{N,\rm per}$ with Dirichlet condition $U_N(x,y)=V_N(x,y)$ when $y=h_N(x)$ is written as
\begin{equation}
\label{eq:UNgU1}
   U_N(x,y) = (\log|\log\delta_N|)\; \gU_N(x,y) \quad\mbox{with}\quad
   \gU_N(x,y) = \frac{1}{N}\gU_1(Nx,Ny)
\end{equation}
where $\gU_1$ is the bounded solution of the Dirichlet problem in $Q_{1,\rm per}$:
\[
   \Delta \gU_1=0\quad\mbox{and}\quad \gU_1(x,y)=y\;\mbox{ for }y=h(x).
\]
The function $\gU_1$ can be considered as the solution of a cell problem associated with a periodic boundary, see Oleinik et al.\ \cite{OSYbook92} and more specifically Nazarov \cite{Nazarov2007}. Formula \eqref{eq:UNgU1} shows how our constructions are completely explicit and why proofs may use nothing more than simple scaling arguments for Sobolev norms when combined with the blow up $\log|\log\delta_N|\to\infty$ as $N\to\infty$.
\medskip

The proof of Theorem \ref{T:NotH32} is performed in sec.\ \ref{S:Conc} by gluing together scaled versions of $\widetilde Q_{N}$ for an increasing sequence $N\sj$, $j\in\N$, and considering in each of them a scaled and truncated version of $U_{N\sj}-V_{N\sj}$. The non-closed range property is then a consequence of our blow up estimates and the closed graph theorem.

\section{Subdomains of the unit disk (circular saws)}
\label{S:disk}
\subsection{The domains $\Omega_{N}$}
\label{SS:OmegaN}
We define a subdomain $\Omega_{1}$ of the unit disk $\Omega_{0}=\D=B_{1}(0)$ in the complex plane described by a graph in polar coordinates
\begin{equation} 
\label{E:Om}
 \Omega_{1} = \{z\in\C\mid \; z=r\re^{i\theta}; \; 0\le r<\rho(\theta)\}\,,
\end{equation}
where $\rho$ is a $2\pi$-periodic Lipschitz continuous and piecewise $\sC^2$ function  satisfying 
$$
  R_{0}\le\rho(\theta)\le1 \quad\mbox{ with some }\; R_{0}>0.
$$
The function $\rho$ can be chosen arbitrarily, the results remain true in any case. It may be piecewise smooth or even analytic as in the specific example described below. The only condition we impose is that {\em it is not constant}, i.e. $\Omega_{1}$ is not a disk centered at the origin.

Let now $N\in\N$, $N\ge1$. We define 
\begin{equation} 
\label{E:OmN}
 \Omega_{N} = \{z\in\C\mid z^{N}\in\Omega_{1}\}\,.
\end{equation}
It is clear that $\Omega_{N}$ has the description by a graph as
$r<\rho_{N}(\theta)$ with $\rho_{N}(\theta)=\rho(N\theta)^{\frac1N}$.

Since the boundaries $\Gamma_{N}$ are described by the equations $r=\rho_{N}(\theta)$, we see from $R_{0}^{\frac1N}\le\rho_{N}(\theta)\le1$ that they converge (in Hausdorff distance) to the unit circle, and from 
$\rho_{N}'(\theta)=\rho(N\theta)^{\frac1N-1}\rho'(N\theta)$ that they are uniformly Lipschitz,  so \eqref{eq:ON} holds. 

\begin{example}
 \label{Ex:Omega}
 For a simple explicit example, we define the domain $\Omega_{1}$ as a conformal (rational) image of the unit disk 
 \begin{equation}
\label{E:Omex}
\Omega^{\Ex}_{1} = \{\omega(z)\mid z\in\D\} \qquad\mbox{ with }\;\;
  \omega(z) = \frac{z}{2-z}\,.
\end{equation}
Then $\Omega^{\Ex}_{N}$ as defined by \eqref{E:OmN} is the image of $\D$ under the conformal mapping $\omega_{N}$ with
$$
\omega_{N}(z) = \omega(z^{N})^{\frac1N} = \frac{z}{(2-z^{N})^{\frac1N}}\,.
$$
It is not hard to see that $\Omega^{\Ex}_{1}$ is the disk $B_{\frac23}(\frac13)$, which implies that we can take $R_{0}=\tfrac13$. 
One can also find the function $\rho$ describing the boundary according to \eqref{E:Om}.
$$
 \rho^{\Ex}(\theta) = \tfrac13(\cos\theta+\sqrt{\cos^{2}\theta+3})\,.
$$
In Figure~\ref{F:1} we show the boundaries $\Gamma^{\Ex}_{N}$ for $N\in\{1,3,9,27\}$ inside the unit circle.
\begin{figure}[h]
\centering
\includegraphics[width=0.6\textwidth]{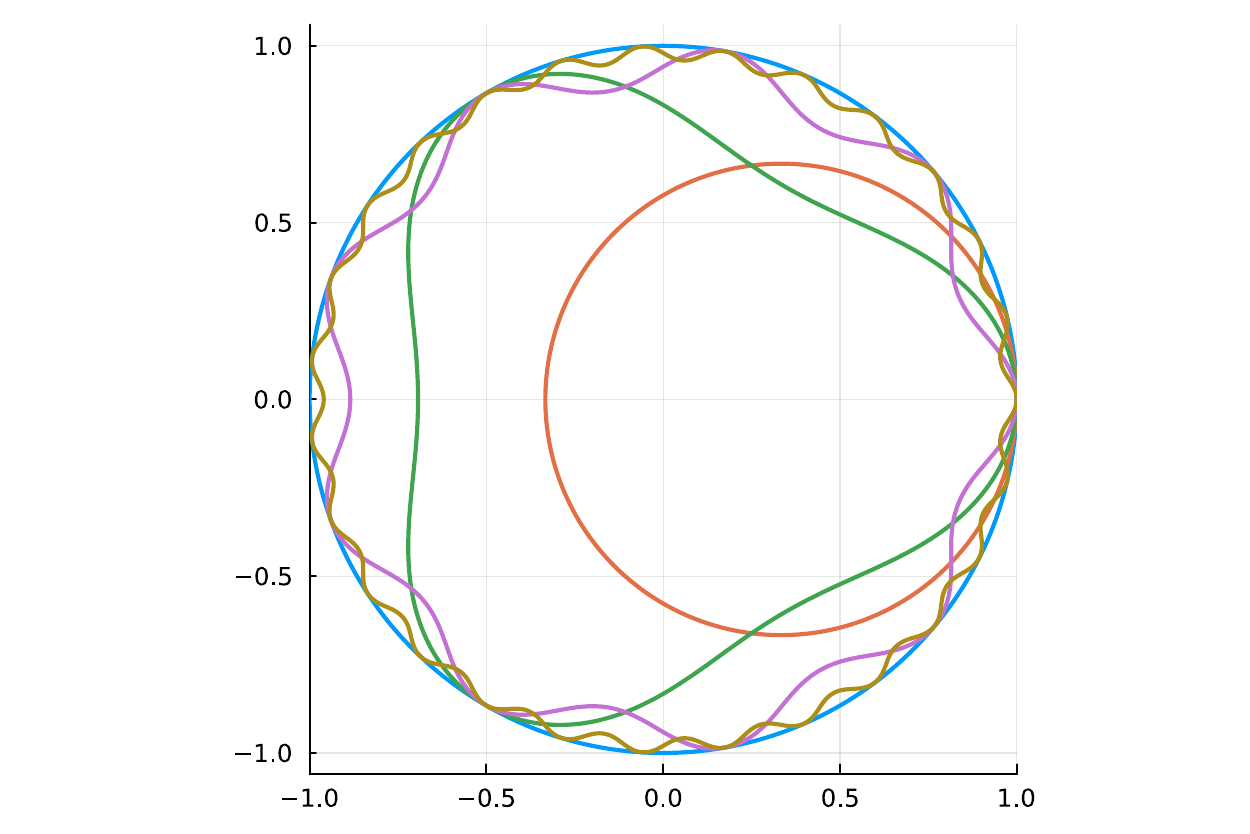}
\caption{Example suite $\Gamma^{\Ex}_{N}$}
\label{F:1}
\end{figure}
For further reference, we note that the Riemann mapping from $\Omega^{\Ex}_{N}$ to $\D$  that maps $0$ to $0$ is given by
$$
 \omega_{N}^{-1}(z) = \frac{2^{\frac1N}z}{(1+z^{N})^{\frac1N}}\,.
$$

\end{example}

\subsection{The right hand side}
\label{SS:rhs}
The sequences $f_{N}$ and $v_{N}$ are based on a function $a\in H^{\frac12}(\R_{+})$ that is unbounded at $0$. We take 
$$
  a(t) = \chi(t)  \log|\log t| \quad\mbox{ with }\chi\in C_{0}^{\infty}([0,\infty)),\; \operatorname{supp}\chi\subset[0,1),\:
  \chi=1 \mbox{ on } [0,\tfrac12]\,.  
$$  
For $\delta>0$ we define
$$
   a_{\delta}(t) = \begin{cases}a(\delta) &\mbox{ for } 0<t<\delta,\\
                              a(t) &\mbox{ for } t\ge\delta.
                  \end{cases}
$$
We can show  that $a_{\delta}\in H^{\frac12}(\R_{+})$ with support contained in $[0,1)$ and 
$$
  \lim_{\delta\to0} \DNorm{a_{\delta}-a}{H^{\frac12}(\R_{+})} = 0\,.
$$
(The easiest way to see this is to write $a(t)$ as trace on $\R_{+}$ of the radial function $a(r)$ on $\R^{2}$ and to prove the corresponding statements for the $H^{1}$-norm of the latter.)

Then we set
$$
  A_{\delta}(t)=\int_{0}^{t} a_{\delta}(s)\,ds\,,
$$
so that $A_{\delta}(t)$ is constant for $t\ge1$, and for $\delta\in(0,\tfrac12)$, $A_{\delta}$ is linear on $[0,\delta]$:
$$
  A_{\delta}(t) = a(\delta)\,t = (\log|\log\delta|)\;t \quad\mbox{ for }\quad 0\le t\le \delta.
$$
 
Then, recalling $R_{0}$ from Section~\ref{SS:OmegaN}: $0<R_{0}\le\min\{\rho(\theta),\theta\in[0,2\pi]\}<1$
, we define for any $N\ge1$ {\em the radial function} $v_N$ given for $z\in\D$ by:
\begin{equation}
\label{E:vN}
 v_{N}(z)=A_{\delta_{N}}(|\log r|)\qquad\mbox{ with }\;
 \delta_{N}=\tfrac1N|\log R_{0}|.
\end{equation} 

Considering that $v_{N}$ is constant in the neighborhood $B_{\frac1\re}(0)$ of the origin and that the mapping 
$$
  \xi^{-1}: z=r\re^{i\theta}\mapsto(x,y)=(\theta,|\log r|)\,:\quad
  \D\setminus\overline{B_{\frac1\re}(0)} \to 
  \T\times(0,1)
$$
which transforms $v_{N}(z)$ to $A_{\delta_{N}}(y)$
is a diffeomorphism, the following is immediately clear:

\begin{lemma}
\label{L:veps}\ 
\\ {\rm (i)} $v_{N}\in H^{\frac32}(\D)$ and 
$\DNorm{v_{N}}{H^{\frac32}(\D)}$ remains bounded as $N\to\infty$.
\\ {\rm (ii)} If $N\ge 2|\log R_{0}|$, (i.e.\ $\delta_N\le\frac12$) then
\begin{equation}
\label{E:vNtrace}
  v_{N}(z) = (\log|\log\delta_N|)\; |\log r| 
  \quad\mbox{ for } \; z\in\Gamma_{N}\,.
\end{equation}
\end{lemma}
  
Finally we set $f_{N}=\Delta v_{N}$ and obtain
\begin{lemma}
\label{L:fN}\ 
\\ {\rm (i)} $\DNorm{f_{N}}{H^{-\frac12}(\D)}$ remains bounded as $N\to\infty$.
\\ {\rm (ii)} $f_{N}$ is radial and is given by 
$f_N(z) = r^{-2} a_{\delta_N}'(|\log r|)$. 
In particular, $f_N$ vanishes for 
$r\in[0, \frac1\re] \cup [\re^{-\delta_N},1]$.
\end{lemma}

\begin{remark}
 \label{R:1-r}
 An entirely equivalent choice for the definition of $v_{N}$ would have been to replace the argument $|\log r|$ in \eqref{E:vN} by $1-r$ and $\delta_{N}$ by $1-\re^{-\delta_{N}}=1-R_{0}^{\frac1N}$.
\end{remark}
 
\subsection{The solution}\label{SS:sol}
As a final tool, we use the solution $\gu_1\in H^{1}(\Omega_1)$ of the Dirichlet problem in the ``generating domain'' $\Omega_1$
\begin{equation}
\label{E:Dir}
  \Delta \gu_1 =0 \quad\mbox{ in } \Omega_1\,,\qquad
  \gu_1(z) = |\log r| \quad\mbox{ on }\Gamma_1 = \partial\Omega_1\,.
\end{equation}
This is the ``corrector'' for obtaining the Green function on $\Omega_1$ from the whole-space fundamental solution: With $r=|z|$
$$
  G(0,z) = -\tfrac1{2\pi}(\log r + \gu_1(z))
$$
is the Green function with pole at the origin.
The function $\gu_1$ is positive harmonic in $\Omega_1$.

\begin{lemma}
\label{L:gu1H32}
The function $\gu_1$ has a finite $H^{\frac32}(\Omega_1)$ norm.
\end{lemma}

Under our assumption that $\rho$ is piecewise $\sC^2$, it is known from the analysis of corner singularities going back to \cite{Kondratev67,Grisvard85,DBook} that there exists even an $\epsilon>0$ such that $\gu_1\in H^{\frac32+\epsilon}(\Omega_1)$. 

\begin{remark}
From the fact that $\Omega_1$ is a Lipschitz domain it is known to follow that the solution $\gu_1$ of the Dirichlet problem with $H^{1}(\Gamma_{1})$ boundary data belongs to $H^{\frac32}(\Omega_1)$, see \cite{JerisonKenig1995}. The standard proof of this  regularity result uses, however, arguments from harmonic analysis that, in this paper, we want to stay clear of. This is the only reason why we are making the stronger assumption that $\rho$ is piecewise $\sC^2$. 
\end{remark}

\smallskip

For the domain $\Omega^{\Ex}_{1}$ of Example~\ref{Ex:Omega}, we can write down $\gu_1$ explicitly: 
Since $\omega^{-1}$ is the Riemann mapping function from $\Omega^{\Ex}_{1}$ to $\D$, the function 
$$
  G^{\Ex}(0,z) = -\tfrac1{2\pi}\log(|\omega^{-1}(z)|)
$$
is the Green function of $\Omega^{\Ex}_1$. Hence
\begin{equation}
\label{E:uNEx}
 \gu^{\Ex}_1(z) = \log(|\omega^{-1}(z)|) - \log r = \log2 - \log|1+z|.
\end{equation}
\qed

Finally, we define
\begin{equation}
\label{eq:guN}
   \gu_{N}(z)=\tfrac1N\gu_1(z^{N})
\end{equation} 
and
\begin{equation}
\label{E:uNdef}
 u_{N} = a(\delta_{N}) \,\gu_{N}\,, \quad\mbox{ i.e. }\quad
 u_{N}(z) = \tfrac1N(\log|\log\delta_N|)\, \gu_1(z^{N})\,.
\end{equation}
The function $\gu_{N}$ is defined  on $\Omega_{N}$ (see \eqref{E:OmN}) and harmonic. On $\Gamma_{N}$ its trace satisfies
$$
  \gu_{N}(z)= \tfrac1N\gu_1(z^{N}) = \tfrac1N(-\log|z^{N}|) = -\log |z|\,. 
$$
Thus $\gu_{N}$ is the corrector function for the Green function of $\Omega_{N}$. 

Comparing with \eqref{E:vNtrace}, we see that $u_{N}=v_{N}$ on $\Gamma_{N}$, hence $u_{N}$ is the solution of the Dirichlet problem \eqref{E:DirN}. 
Our aim is, relying on Lemma \ref{L:gu1H32}, to show that the $H^{\frac32}$ norm of $\gu_{N}$ is bounded from below by a positive constant independent of $N$, which implies that the $H^{\frac32}$ norm of $u_{N}$ will tend to infinity as $N\to\infty$.
\smallskip

\begin{figure}[h]
\centering
\includegraphics[width=0.3\textwidth]{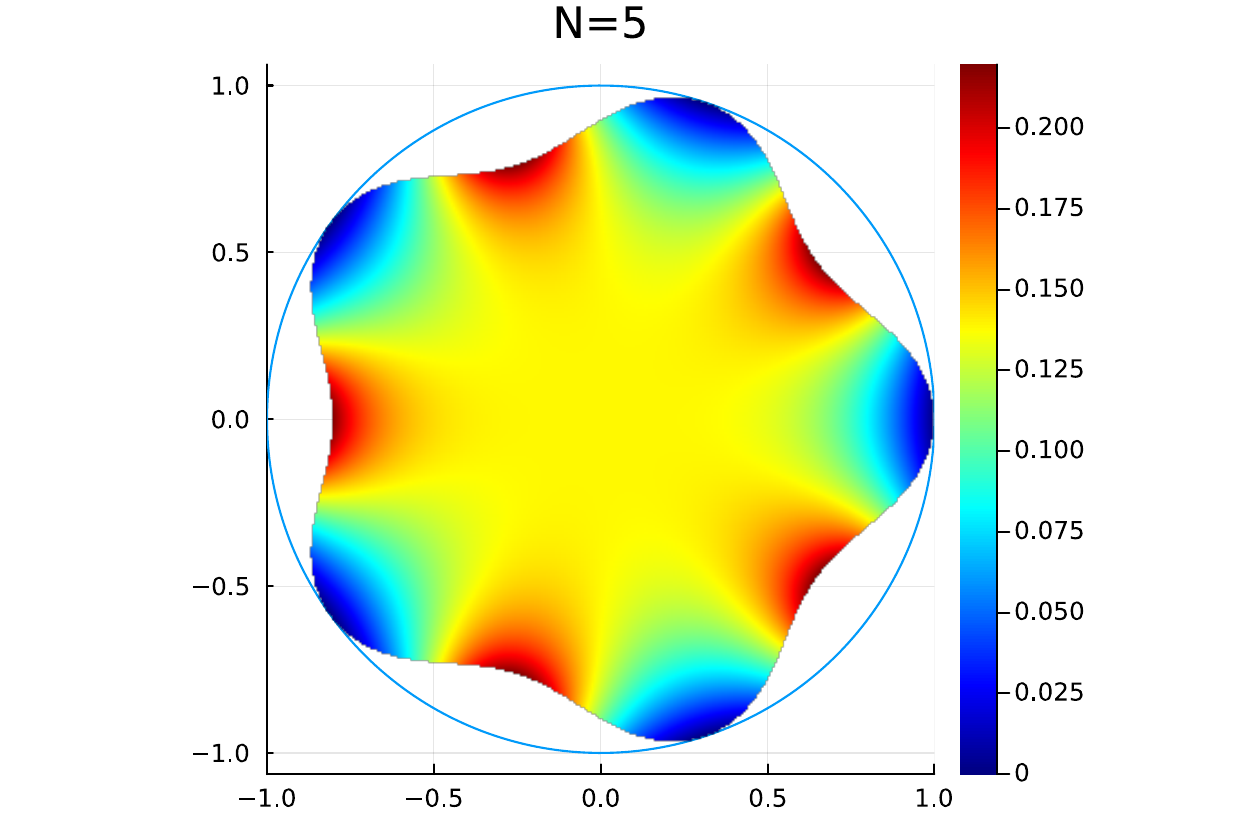}\ 
\includegraphics[width=0.3\textwidth]{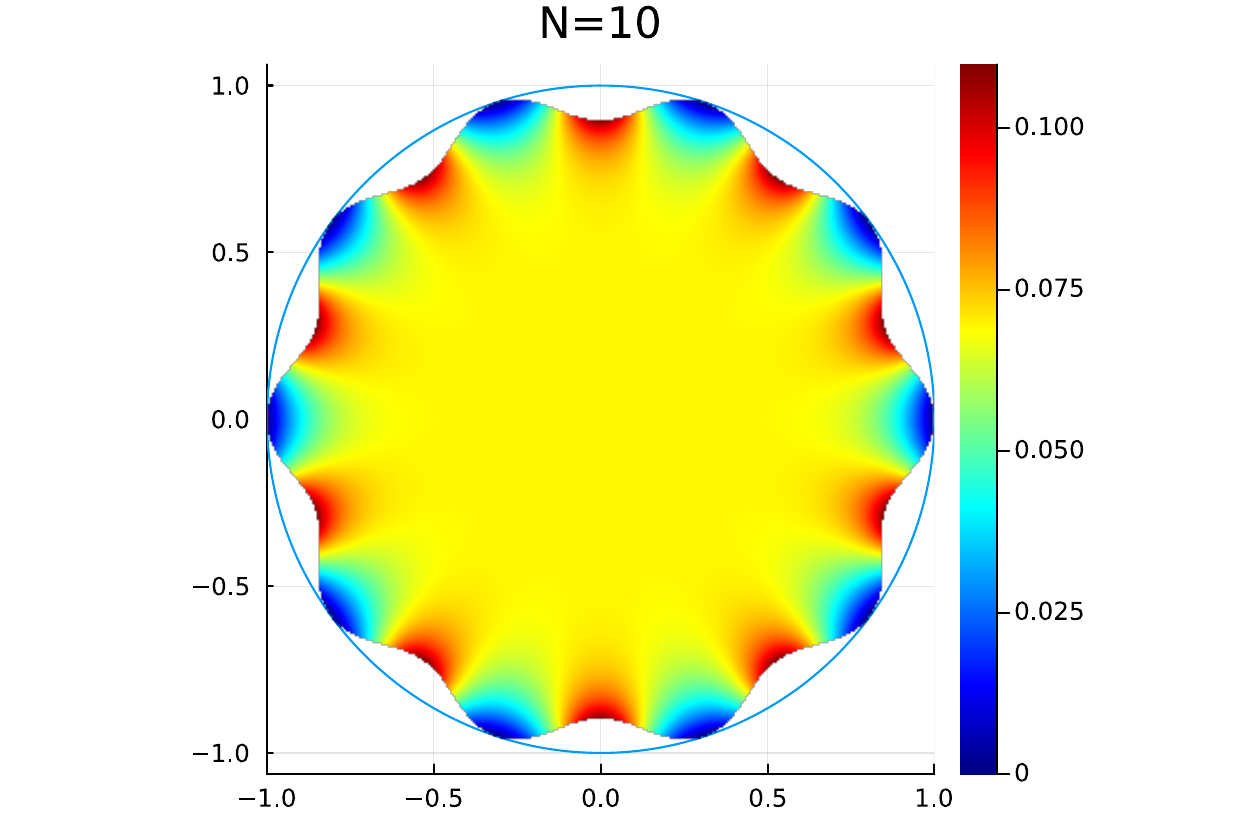}\ 
\includegraphics[width=0.3\textwidth]{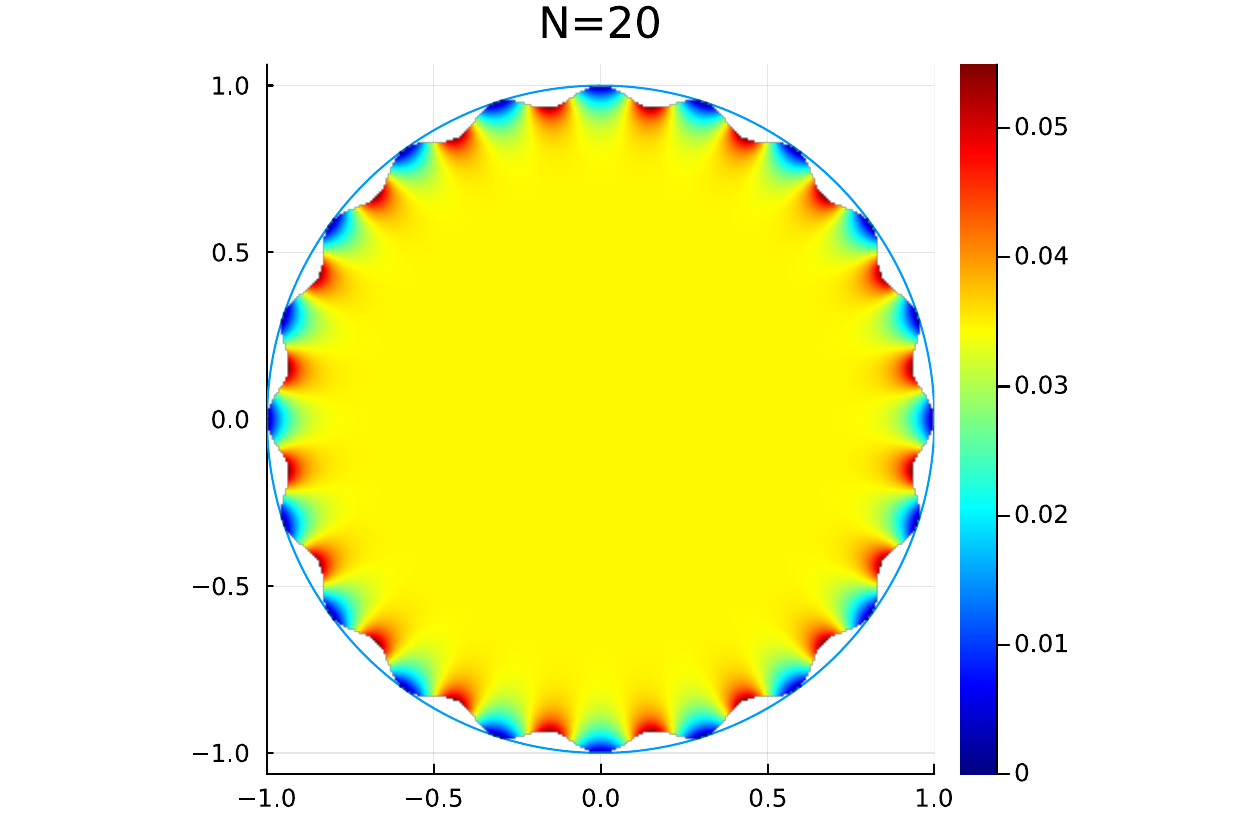}
\caption{Example solution $\gu^{\Ex}_{N}$}
\label{F:2}
\end{figure}

In Figure~\ref{F:2}, we plot $\gu^{\Ex}_{N}=\tfrac1N\gu^{\Ex}_1(z^{N})$ for several values of $N$. The boundary layer structure is clearly visible.

\section{Local norm estimates in circular saws}\label{S:est}
In order to show that the $H^{\frac32}(\Omega_N)$ norm of $\gu_{N}$ \eqref{eq:guN} is bounded from below by a positive constant independent of $N$, and hence to prove Proposition \ref{P:a}, it suffices to prove this on appropriate subdomains of $\Omega_N$, possibly depending on $N$. We will see in this section that we can consider annular regions of $\Omega_N$ close to the boundary $\Gamma_N$ with radii $R^{\frac1N}$ and $(R')^{\frac1N}$ for $0<R<R'\le R_0$.

\subsection{Change of variables}\label{SS:chvar}
Cartesian coordinates being easier to manipulate, we employ the conformal mapping (already mentioned in the introduction)
\begin{equation}
\label{E:}
 \xi: (x,y)\mapsto r\re^{i\theta} \quad\mbox{ with } r=\re^{-y}, \theta=x\,,
 \qquad\mbox{i.\ e.}\quad \xi(x,y)=\re^{i(x+iy)}\,.
\end{equation}
Under this mapping, for the domain $\Omega_{N}$ as described in \eqref{E:OmN},
$\Omega_{N}\setminus[0,1)$ is the image of the semi-infinite strip
\begin{equation}
\label{E:qN}
 Q_{N} = \{(x,y)\mid 0< x<2\pi ; \;\;y>h_{N}(x)\}\,, 
\end{equation}
where
\begin{equation}
\label{eq:hN}
 h_{N}(x)=\tfrac1N h(Nx)  \quad\mbox{with}\quad h(x) = - \log\rho(x)\,,
\end{equation}
and $\Omega_{N}$ is the image of $Q_{N,\rm per}$, which is $Q_{N}$ with the boundary at $x=0$ and $x=2\pi$ added and identified.

\begin{figure}[h]
\centering
\includegraphics[width=0.48\textwidth]{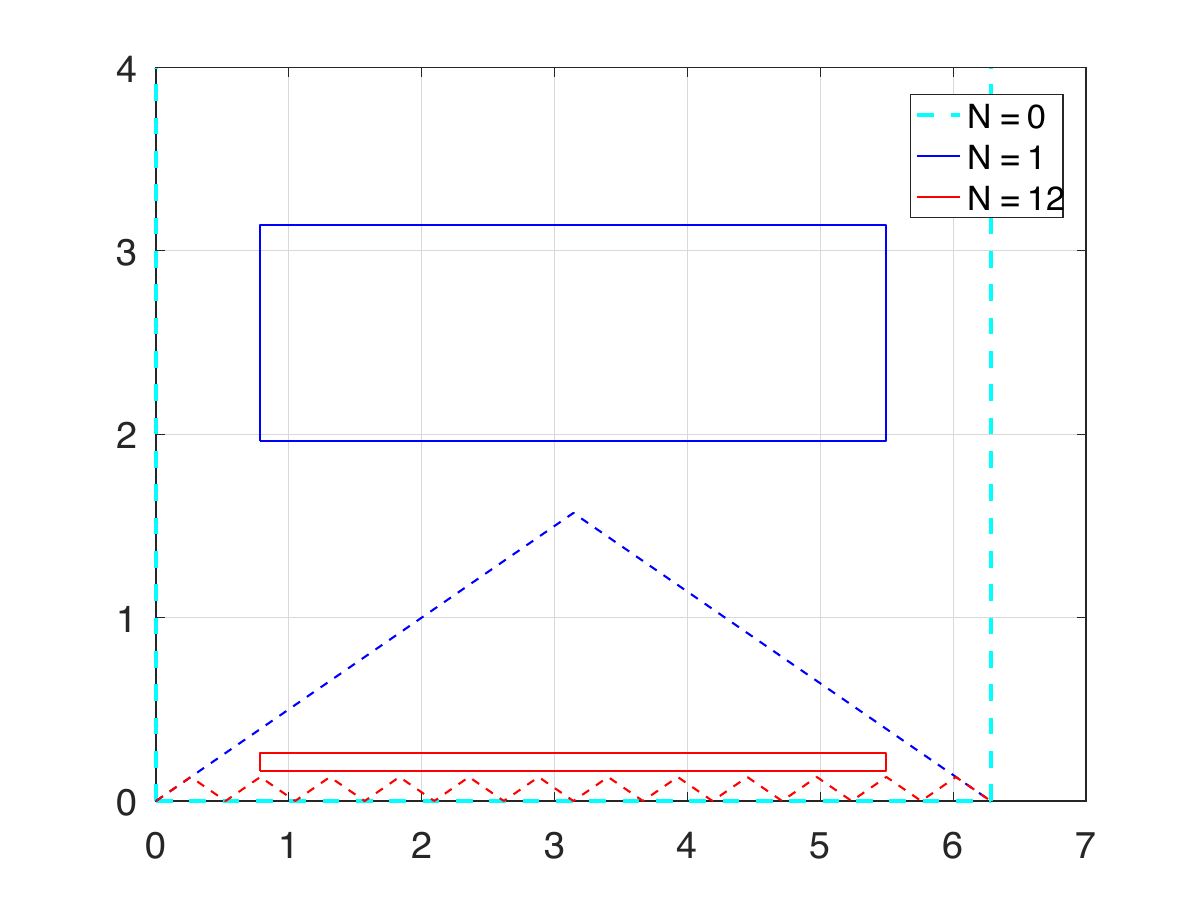}\quad
\includegraphics[width=0.48\textwidth]{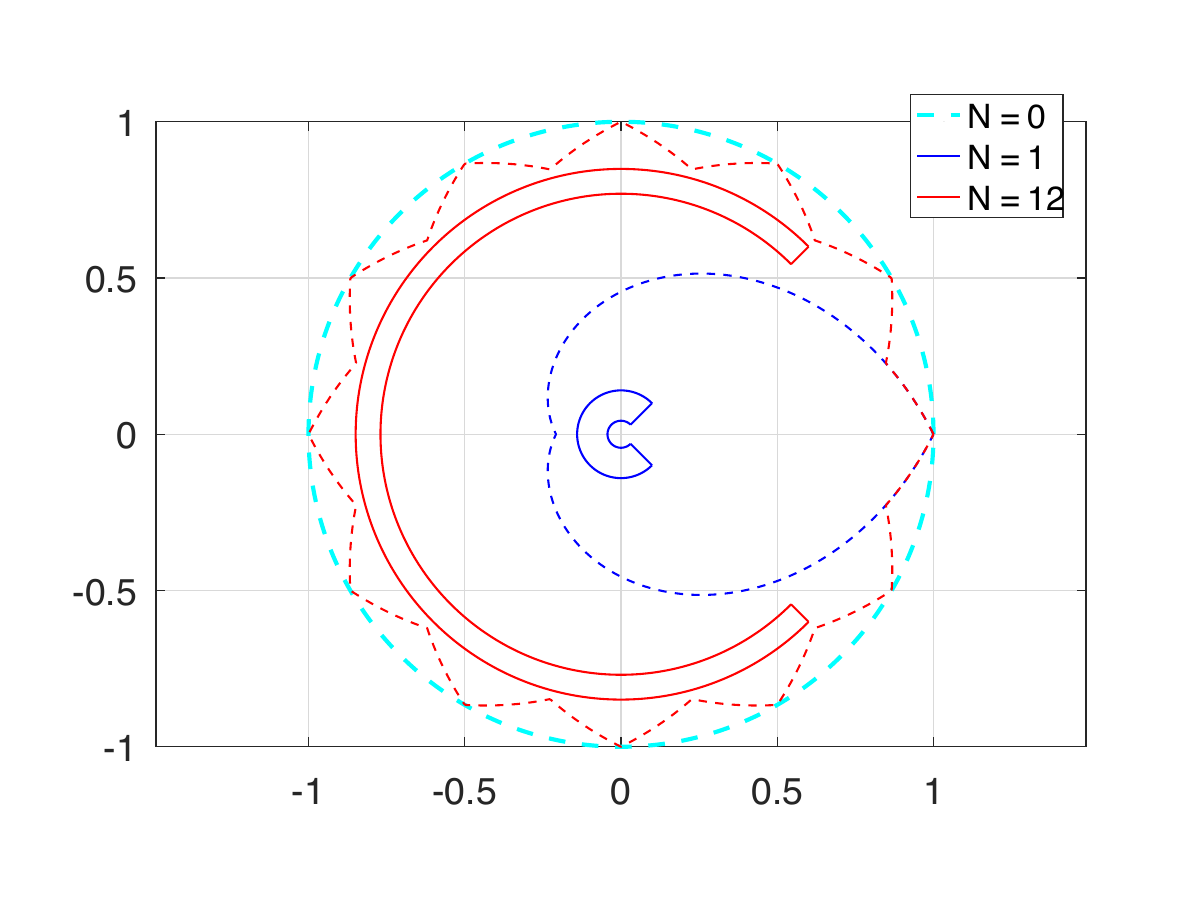}
\caption{\small Domains $Q^{a,b,\alpha,\beta}_{N}$ (left) and $\Omega^{T,T',R,R'}_{N}$ (right) plotted with solid lines for $N=1$ and $12$ when the motif $h$ is the the half of the hat function defined in \eqref{eq:hat}, and the values $a=\frac\pi4$, $b=2\pi-a$, $\alpha=\frac58\pi$, $\beta=\pi$ are chosen. Corresponding domains $Q_N$ and $\Omega_N$ are plotted with dotted lines.}
\label{F:hatrect}
\end{figure}

For the ease of later calculations we introduce the cartesian rectangular subdomain of $Q_N$
\begin{equation}
\label{E:QN}
 Q^{a,b,\alpha,\beta}_{N} = 
   (a,b)\times(\tfrac{\alpha}{N},\tfrac{\beta}{N})
   \quad
   \mbox{ with }\quad
   0\le a<b<2\pi,\quad \max h<\alpha < \beta \,,
\end{equation}
which is in one by one correspondence with
the rectangle in polar coordinates
\begin{equation}
\label{E:OmTR}
 \Omega^{T,T',R,R'}_{N} = 
  \{r\re^{i\theta}\mid T<\theta<T', \quad R^{\frac1N}<r<(R')^{\frac1N}\} \,,
\end{equation}
with $a=T$, $b=T'$, $R=\re^{-\beta}>0$, and $R'=\re^{-\alpha}\in (R,R_0)$.
See Figure~\ref{F:hatrect} for an illustration.

For the simplicity of norm estimates, it will be useful to be able to switch between both systems of coordinates. 

Since the mapping $\xi$ is $C^{\infty}$ everywhere and its inverse is $C^{\infty}$ except at the origin, and thus $\xi^{-1}$ has bounded derivatives of any order on $\D\setminus B^{R}(0)$ for $R>0$, we immediately get equivalence of Sobolev norms in both systems of coordinates:
 \begin{equation}
\label{E:cart-polar}
  c \DNorm{u}{H^{s}(\Omega^{T,T',R,R'}_{N})} \le
   \DNorm{u\circ\xi}{H^{s}(Q^{a,b,\alpha,\beta}_{N})} \le 
  C \DNorm{u}{H^{s}(\Omega^{T,T',R,R'}_{N})}\,.
\end{equation}
The constants $c$ and $C$ may depend on $R$ and $s$, but not on $T,T',R'$, and $u$.

\subsection{Scaling}\label{SS:scal}
Recall the function $\gu_{N}$ \eqref{eq:guN} from Section~\ref{SS:sol}. The function 
$$
 \gU_{N} = \gu_{N}\circ\xi \,: \quad \gU_{N}(x,y) = \gu_{N}(\re^{i(x+iy)})
$$
is harmonic, $2\pi$-periodic in $x$, and defined for $y>h_{N}(x)$, where
$h_{N}(x)=\tfrac1N h(Nx)$, $h(x)=-\log\rho(x)$.
It satisfies the scaling relation
\begin{equation}
\label{E:scale}
 \gU_{N}(x,y)=\tfrac1N \gU_{1}(Nx,Ny)\,.
\end{equation}

We use the Sobolev seminorms on a domain $\cO\subset\R^{2}$ defined as usual by 
$$
  \Normc{f}{H^{s}(\cO)}2 = \sum_{|j|=m}\int_{\cO}|D^{j}f(z)|^{2}dz
$$
for $s=m\in\N$ and
$$
  \Normc{f}{H^{s}(\cO)}2 = \sum_{|j|=m}\Normc{D^{j}f}{H^{\sigma}(\cO)}2
$$
for $s=m+\sigma$, $m\in\N$, $0<\sigma<1$, where the Slobodetsky seminorm is defined as
$$
 \Normc{f}{H^{\sigma}(\cO)}2 = 
  \int_{\cO}\int_{\cO}|f(z)-f(z')|^{2}|z-z'|^{-2-2s} \,dz\,dz'\,,
$$
and then the norms by
$$
  \DNormc{f}{H^{s}(\cO)}2 = \sum_{\ell\le m}\Normc{f}{H^{\ell}(\cO)}2
  \quad\text{ for }s=m\in\N\,,
$$
$$
  \DNormc{f}{H^{s}(\cO)}2 = \DNormc{f}{H^{m}(\cO)}2 + \Normc{f}{H^{\sigma}(\cO)}2
  \quad\text{ for }s=m+\sigma\,.
$$
In addition, we need the norm of the periodic Sobolev space 
$H^{s}_{\rm per}((0,2\pi)\times (\alpha,\beta))$ for functions that are $2\pi$-periodic in the first variable. This is obtained by replacing the euclidean distance $|z-z'|$ in the Slobodetsky integral by the distance $|z-z'|_{\rm per}$ on the flat torus $\R/(2\pi\Z)\times (\alpha,\beta))$:
$$
 |(x,y)-(x',y')|_{\rm per}^{2}= |x-x'|_{\rm per}^{2}+ (y-y')^{2},
 \quad |x-x'|_{\rm per} = \min\{|x-x'|,2\pi-|x-x'|\}\,.
$$
\begin{equation}
\label{E:Hsper}
 \Normc{f}{H^{\sigma}_{\rm per}((0,2\pi)\times (\alpha,\beta))}2 = 
  \int_{0}^{2\pi}\!\!\int_{\alpha}^{\beta}\int_{0}^{2\pi}\!\!\int_{\alpha}^{\beta}
  \frac{|f(x,y)-f(x',y')|^{2}}{|(x,y)-(x',y')|_{\rm per}^{2+2\sigma}}
  \,dy\,dx\,dy'\,dx'
\end{equation}
and 
$\DNormc{f}{H^{s}_{\rm per}(\cO)}2 = \DNormc{f}{H^{m}(\cO)}2 + \Normc{f}{H^{\sigma}_{\rm per}(\cO)}2$. It is straightforward to verify that for a function $f$ that is $2\pi$-periodic in the first variable, we have the following inclusion between periodic and non-periodic seminorms:
\begin{equation}
\label{E:per-nonper}
 \Norm{f}{H^{\sigma}((0,2\pi)\times (\alpha,\beta))}
 \le \Norm{f}{H^{\sigma}_{\rm per}((0,2\pi)\times (\alpha,\beta))}
  \le \Norm{f}{H^{\sigma}((0,4\pi)\times (\alpha,\beta))}\,.
\end{equation}

The transformation $\xi$ to cartesian variables allows us to benefit from the  fact that Sobolev seminorms are homogeneous with respect to change of scale (dilations) in the independent variable.

\begin{lemma}
 \label{L:scale}
Let $U_1$ be smooth in $\T\times[\alpha,\beta]$ and $2\pi$-periodic in the first variable. Let $U_N$ satisfy the scaling relation \eqref{E:scale} $U_{N}(x,y)=\tfrac1N U_{1}(Nx,Ny)$. Then
 for any $s\ge0$
 \begin{equation}
\label{E:scaleHs}
 \Norm{U_{N}}{H^{s}(Q^{a,b,\alpha,\beta}_{N})} 
  = N^{s-2}\,\Norm{U_{1}}{H^{s}(Q^{Na,Nb,\alpha,\beta}_{1})} .
\end{equation}
\end{lemma}

\begin{proof}
 Note that $U_{N}$ is $C^{\infty}$ on this rectangle, so we can simply change variables in the integrals defining these seminorms, using the scaling law \eqref{E:scale}. For integer $s$ and $|j|=s$ we find
$$
 \begin{aligned}
\int_{a}^{b}\int_{\frac\alpha{N}}^{\frac\beta N} |D^{j}U_{N}(x,y)|^{2} \,dy\,dx &=
 \int_{a}^{b}\int_{\frac\alpha{N}}^{\frac\beta N} |N^{s-1}(D^{j}U_{1})(Nx,Ny)|^{2} \,dy\,dx\\
 &= N^{2s-4} \int_{Na}^{Nb}\int_{\alpha}^{\beta} |D^{j}U_{1}(X,Y)|^{2} \,dY\,dX
 \end{aligned}
$$
Similarly, for $0<s<1$
$$
 \begin{aligned}
\int_{Q^{a,b,\alpha,\beta}_{N}}\int_{Q^{a,b,\alpha,\beta}_{N}}&
  \frac{|U_{N}(x,y)-U_{N}(x',y')|^{2}}{|(x,y)-(x',y')|^{2+2s}}\,dx\,dx'\,dy\,dy'
  \\ &=
  N^{2s}
 \int_{Q^{a,b,\alpha,\beta}_{N}}\int_{Q^{a,b,\alpha,\beta}_{N}}
  \frac{|U_{1}(Nx,Ny)-U_{1}(Nx',Ny')|^{2}}{|(Nx,Ny)-(Nx',Ny')|^{2+2s}} \,dx\,dx'\,dy\,dy'\\
 &= N^{2s-4} \int_{Q^{Na,Nb,\alpha,\beta}_{1}}\int_{Q^{Na,Nb,\alpha,\beta}_{1}}
  \frac{|U_{1}(X,Y)-U_{1}(X',Y')|^{2}}{|(X,Y)-(X',Y')|^{2+2s}} \,dX\,dX'\,dY\,dY'
 \end{aligned}
$$
\end{proof}

Next we use the $2\pi$-periodicity of $U_{1}$ in $x$
to estimate the seminorm on 
$Q^{Na,Nb,\alpha,\beta}_{1}=(Na,Nb)\times(\alpha,\beta)$.
We are going to prove:

\begin{lemma}
\label{L:U1per}
For $U_1$ as in Lemma \ref{L:scale}, we have
\begin{equation}
\label{E:HsN-per}
c \,N \Normc{U_{1}}{H^{s}_{\rm per}((0,2\pi)\times (\alpha,\beta))}2
\le \Normc{U_{1}}{H^{s}(Q^{Na,Nb,\alpha,\beta}_{1})}2
\le
C \,N \Normc{U_{1}}{H^{s}_{\rm per}((0,2\pi)\times (\alpha,\beta))}2
\end{equation}
with constants $c$, $C$ independent of $N$.
\end{lemma}

\begin{proof}
Let $k\in\N$ be such that $2\pi k\le N(b-a)<2\pi(k+1)$. 
Then 
$$
  k\int_{0}^{2\pi}|D^{j}U_{1}(X,Y)|^{2}\,dX \le
  \int_{Na}^{Nb}|D^{j}U_{1}(X,Y)|^{2}\,dX \le
 (k+1)\int_{0}^{2\pi}|D^{j}U_{1}(X,Y)|^{2}\,dX .
$$
Hence for integer $s=m$ we find the sharp asymptotics
\begin{equation}
\label{E:quotper}
\lim_{N\to\infty}
  \frac{\Normc{U_{1}}{H^{m}(Q^{Na,Nb,\alpha,\beta}_{1})}2}%
  {\frac{b-a}{2\pi}N\Normc{U_{1}}{H^{m}(Q^{0,2\pi,\alpha,\beta}_{1})}2}
  =1\,,
\end{equation}
which proves \eqref{E:HsN-per} for integer $s$.
\smallskip

For non-integer $s$, we have to take into account that the Sobolev-Slobodetsky double integral has a non-periodic integrand. 

To prove the lower bound in \eqref{E:HsN-per}, we assume that $N(b-a)\ge4\pi k$ and observe for all $Y,Y'$ 
$$
 \begin{aligned}
 \int_{0}^{4\pi k}\int_{0}^{4\pi k} &
  \frac{|U_{1}(X,Y)-U_{1}(X',Y')|^{2}}{|(X,Y)-(X',Y')|^{2+2s}} \,dX\,dX'
  \\ &=
  \sum_{\ell=1}^{k}\sum_{\ell'=1}^{k}
  \int_{0}^{4\pi}\int_{0}^{4\pi}
  \frac{|U_{1}(X,Y)-U_{1}(X',Y')|^{2}}{\big((X-X'+4\pi(\ell-\ell'))^{2}+(Y-Y')^{2}\big)^{1+s}} \,dX\,dX'
  \\ &\ge
    k \int_{0}^{4\pi}\int_{0}^{4\pi}
  \frac{|U_{1}(X,Y)-U_{1}(X',Y')|^{2}}{|(X,Y)-(X',Y')|^{2+2s}} \,dX\,dX'\,, \end{aligned}
$$
so that
$$
 \begin{aligned}
 \Normc{U_{1}}{H^{s}(Q^{Na,Nb,\alpha,\beta}_{1})}2 
  &\ge
   \Normc{U_{1}}{H^{s}(Q^{0,4\pi k,\alpha,\beta}_{1})}2 
  \\ &\ge
   k \Normc{U_{1}}{H^{s}(Q^{0,4\pi,\alpha,\beta}_{1})}2 
  \\ &\ge
  k \Normc{U_{1}}{H^{s}_{\rm per}((0,2\pi)\times (\alpha,\beta))}2 \,.
 \end{aligned}
$$
Here we used \eqref{E:per-nonper} for the last inequality. We obtain the lower bound for any $c<\frac{b-a}{4\pi}$.
\smallskip

The upper bound in \eqref{E:HsN-per} is more delicate. We assume $N(b-a)\le2\pi k$ and  split the integration domain $(0,2\pi k)\times(0,2\pi k)$ into squares of size 
$(0,2\pi)\times(0,2\pi)$:
$$
 \int_{0}^{2\pi k}\int_{0}^{2\pi k} 
  \frac{|U_{1}(X,Y)-U_{1}(X',Y')|^{2}}{|(X,Y)-(X',Y')|^{2+2s}} \,dX\,dX'
 =
  \sum_{\ell=1}^{k}\sum_{\ell'=1}^{k} J_{\ell-\ell'}
  = \sum_{m=-k}^{k}(k-|m|)J_{m}
  \,,
$$
with 
$$
J_{m} = \int_{0}^{2\pi}\int_{0}^{2\pi}
  \frac{|U_{1}(X,Y)-U_{1}(X',Y')|^{2}}{\big((X-X'+2\pi m)^{2}+(Y-Y')^{2}\big)^{1+s}} \,dX\,dX'\,.
$$
For $m=-1,0,1$, we have
$$
 (X-X'+2\pi m)^{2}+(Y-Y')^{2} \ge |(X,Y)-(X',Y')|_{\rm per}^{2}\,,
$$
see \eqref{E:Hsper}. Hence after integration over $Y,Y'$, these terms give a contribution bounded by
\begin{equation}
\label{E:|m|le1}
 (3k-2) \Normc{U_1}{H^{s}_{\rm per}((0,2\pi)\times (\alpha,\beta))}2\,.
\end{equation}
For $|m|\ge2$, we can write
$$
  (X-X'+2\pi m)^{2}+(Y-Y')^{2} 
  = \big( |(X,X')-(Y,Y')|^{2}\big)\big(
    1+ \frac{(X-X'+2\pi m)^{2} - (X-X')^{2}}{|(X,X')-(Y,Y')|^{2}}
    \big)\,.
$$
With $|X-X'|\le2\pi$ and $|Y-Y'|\le \beta-\alpha$, we get 
$$
  (X-X'+2\pi m)^{2}+(Y-Y')^{2} 
  \ge \big( |(X,X')-(Y,Y')|^{2}\big)\big(
    1+ \tfrac{(2\pi)^{2}|m|(|m|-2)}{D}
    \big)\,,
$$
where $D=(2\pi)^{2}+(\beta-\alpha)^{2}$.

Hence after after integration over $Y,Y'$, the terms with $|m|\ge2$ give a contribution bounded by
\begin{equation}
\label{E:|m|ge2}
 2k \big(\sum_{m=2}^{k}(1+\tfrac{(2\pi)^{2}m(m-2)}{D})^{-1-s}\big)
 \Normc{U_1}{H^{s}_{\rm per}((0,2\pi)\times (\alpha,\beta))}2\,.
\end{equation}
The latter sum is bounded independently of $k$ as part of a convergent infinite series.

Together, \eqref{E:|m|le1} and \eqref{E:|m|ge2} give an upper bound 
$C N \Normc{U_1}{H^{s}_{\rm per}((0,2\pi)\times (\alpha,\beta))}2$ with a constant $C$ independent of $N$, and \eqref{E:HsN-per} is proved.
\end{proof}

Altogether, we have proved the following consequence of Lemmas~\ref{L:scale}-\ref{L:U1per}
\begin{lemma}
 \label{L:fullperiod}
 For any $s\ge0$, $a<b$, and $|\log R_{0}|\le\alpha<\beta$, 
 there exist positive constants $c$ and $C$ such that for any 
 $N\ge\tfrac{4\pi}{b-a}$ 
 \begin{equation}
\label{E:estHs}
  c\,\Norm{\gU_{1}}{H^{s}_{\rm per}(Q^{0,2\pi,\alpha,\beta}_{1})} \le
  N^{\frac32-s}
 \Norm{\gU_{N}}{H^{s}(Q^{a,b,\alpha,\beta}_{N})} 
  \le C\,\Norm{\gU_{1}}{H^{s}_{\rm per}(Q^{0,2\pi,\alpha,\beta}_{1})} .
\end{equation} 
Here the lower and upper bounds do not depend on $N$.
\end{lemma}

\subsection{Conclusion}\label{SS:concON}
In order to have a non-zero lower bound in \eqref{E:estHs}, we have to be sure that the $H^{s}$-seminorm of $\gU_{1}$ does not vanish.

\begin{lemma}
\label{L:nonconst}
If $\rho$, hence $h$, are non-constant, then $\Norm{\gU_{1}}{H^{s}_{\rm per}(Q^{0,2\pi,\alpha,\beta}_{1})}$ is not $0$.
\end{lemma}

\begin{proof}
If the $H^{s}$-seminorm of $\gU_{1}$ was $0$, we would find that $\gU_1$ is a polynomial, harmonic and periodic in $x$. This implies that it is of the form $\gU_{1}(x,y)=c+dy$, hence $\gu_1(z)=c-d\log r$. Since $\gu_1$ is regular at the origin, this leaves the only possibility that $\gu_1$ is constant, and in view of its boundary condition \eqref{E:Dir}, this would mean that $\Gamma_1$ is a circle centered at the origin, which we have excluded.
\end{proof}

Before coming back to our annular domains $\Omega^{T,T',R,R'}_{N}$, 
we have to note that the seminorm estimate (relying on \eqref{E:estHs} and Lemma \ref{L:nonconst})
$$
    c\, N^{s-\frac32}\le
 \Norm{\gU_{N}}{H^{s}(Q^{a,b,\alpha,\beta}_{N})} 
  \le C\,N^{s-\frac32} , \qquad \forall N\ge\tfrac{4\pi}{b-a}\,,
$$
together with the corresponding upper bound for $s=0$ implies that one has the same type of estimate for the norms:
$$
    c\, N^{s-\frac32}\le
 \DNorm{\gU_{N}}{H^{s}(Q^{a,b,\alpha,\beta}_{N})} 
  \le C\,N^{s-\frac32} \,.
$$
This one can then be transferred to $\Omega^{T,T',R,R'}_{N}$ via the diffeomorphism $\xi$  \eqref{E:cart-polar}. The result is

\begin{proposition}
 \label{P:estHs}
 Let $\rho$ be non-constant. 
 Let $0\le T< T'\le2\pi$ and $0<  R<R'\le R_{0}$. Then\\
 {\rm (i)} for any $s\ge0$ for which $\gu_1\in H^{s}(\Omega_{1})$, there are positive constants $c$ and $C$ such that for all $N\ge\tfrac{4\pi}{T'-T}$ 
\begin{equation}
   \label{E:normHs}
    c\, N^{s-\frac32}\le
 \DNorm{\gu_{N}}{H^{s}(\Omega^{T,T',R,R'}_{N})} 
  \le C\,N^{s-\frac32} \,.
\end{equation}
{\rm (ii)} In particular, with $u_{N}=a(\delta_{N})\,\gu_{N}$ the solution of the Dirichlet problem \eqref{E:DirN}, see \eqref{E:uNdef}, we have the lower bound for all $N\ge\tfrac{4\pi}{T'-T}$ 
\begin{equation}
\label{E:3/2lower}
 \DNorm{u_{N}}{H^{\frac32}(\Omega_{N})} \ge c\,\log|\log\delta_N| =
 c\,\log(\log(\tfrac N{|\log R_{0}|}))\,.
\end{equation}
\end{proposition}

Thus $\DNorm{u_{N}}{H^{\frac32}(\Omega_{N})} \to \infty$, and, in view of Lemma \ref{L:veps}, the proof of Proposition \ref{P:a} is complete. Hence we have proved Theorem~\ref{T:blowup}.

\section{Localized norm estimates in straight saws}\label{S:Rect}
In this section our aim is twofold: First prove Proposition \ref{P:bc}, second prepare for the proof of Theorem \ref{T:NotH32}. To prove this proposition about periodic and straight saws $\widetilde Q_{N,\rm per}$ and $\widetilde Q_{N}$ we have to truncate functions $W_N=U_N-V_N$ where $U_N=u_N\circ\xi$, $V_N=v_N\circ\xi$, and for this we need global estimates on $u_N$. Such a truncation will also provide elementary bricks for the proof of Theorem \ref{T:NotH32}.

\subsection{Global estimates in $\Omega_N$}\label{SS:global}
\begin{figure}[h]
\centering
\includegraphics[width=0.48\textwidth]{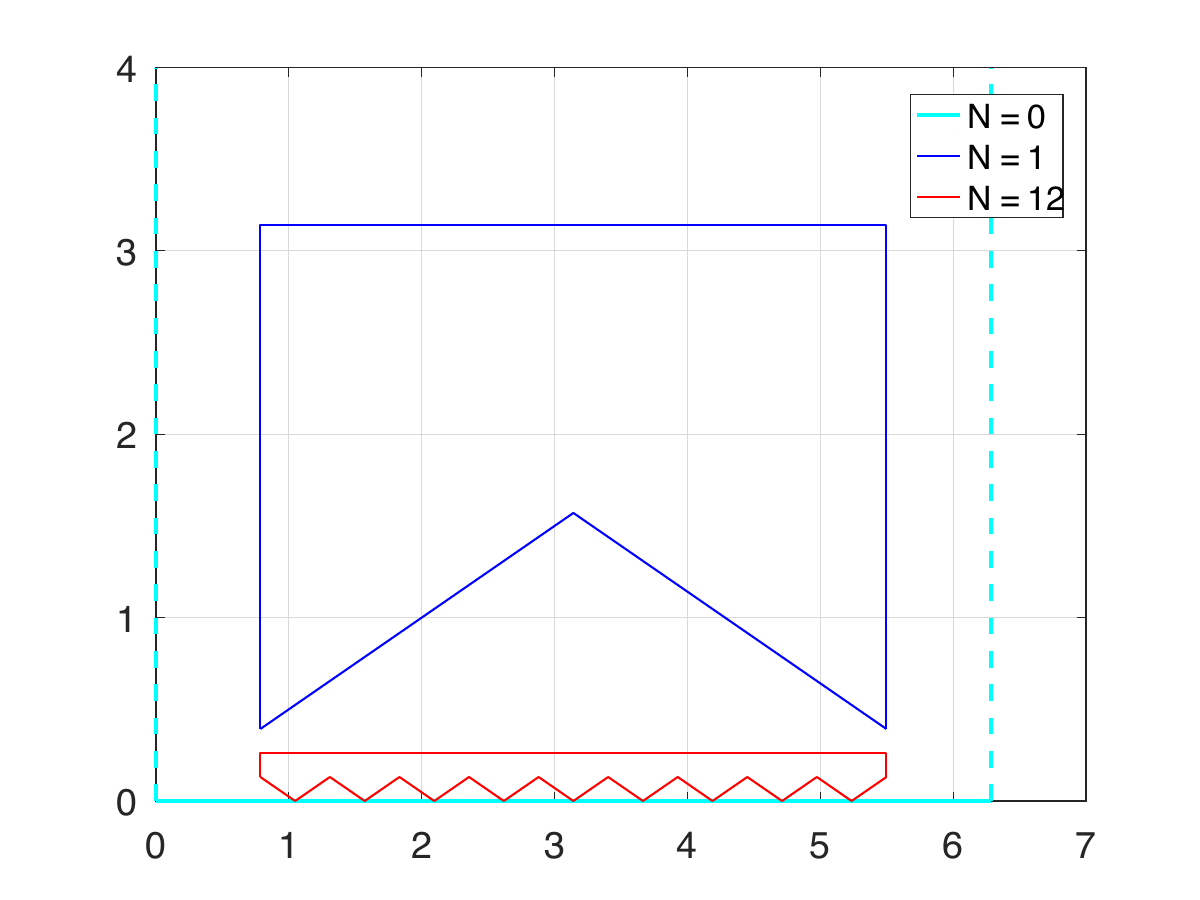}\quad
\includegraphics[width=0.48\textwidth]{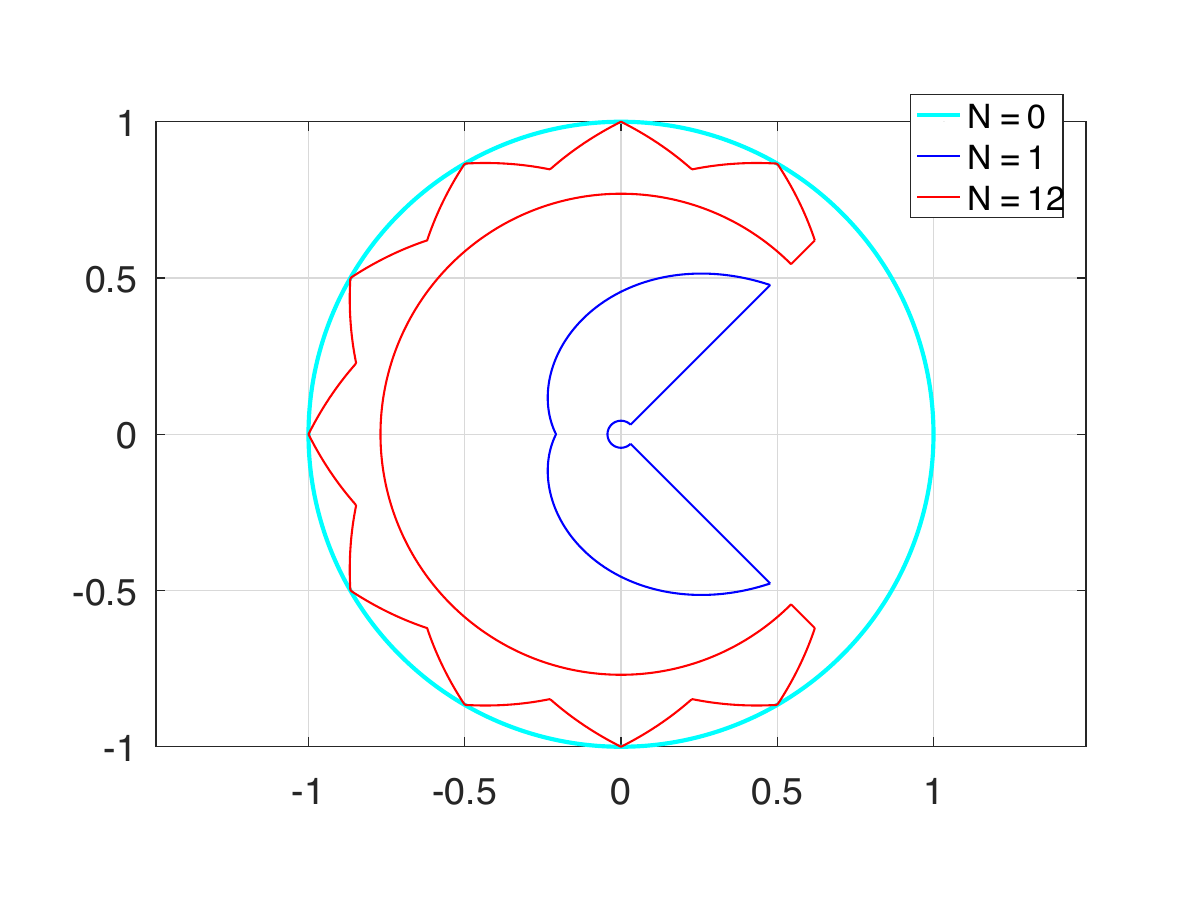}
\caption{\small Domains $Q^{\BL a,b;\,\alpha}_{N}$ (left) and $\Omega^{\BL T,T';\,R}_{N}$ (right) for $N=1$ and $12$ when the motif $h$ is the the half of the hat function defined in \eqref{eq:hat}, and the values $a=\frac\pi4$, $b=2\pi-a$, $\alpha=\pi$ are chosen.}
\label{F:wtrect}
\end{figure}

Now we use layers in $\Omega_N$ up to the boundary which we define in the same spirit as \eqref{E:OmTR}:
\begin{equation}
\label{E:BLOmTR}
  \Omega^{\BL T,T';\, R}_{N} = 
  \{r\re^{i\theta}\mid T<\theta<T', \quad R^{\frac1N}<r<\rho_N(\theta)\} \,,
\end{equation}
for $0\le T<T'\le 2\pi$ and $R<R_0$. The boundaries of these domains contain the part of $\Gamma_N$ corresponding to $\theta\in(T,T')$, while $\Omega^{T,T',R,R'}_{N}$ have no common boundary with $\Omega_N$. On the cartesian side we define, cf \eqref{E:QN}:
\begin{equation}
\label{E:BLQN}
  Q^{\BL a,b;\,\alpha}_{N} = 
  \{(x,y)\mid a<x<b,\quad h_N(x)<y<\tfrac{\alpha}{N}\}
\end{equation}
with $0\le a<b<2\pi$ and $\alpha > \max h$. 
For $a=T$, $b=T'$, and $R=\re^{-\alpha}>0$, the subdomains $\Omega^{\BL T,T';\, R}_{N}$ and $Q^{\BL a,b;\,\alpha}_{N}$ are isomorphic.
See Figure~\ref{F:wtrect} for an illustration choosing the same parameters as in Figure~\ref{F:hatrect}.

With essentially the same proof as in Section~\ref{S:est}, one can prove the same estimates when the domain $\Omega^{T,T',R,R'}_{N}$ is replaced by the layer $\Omega^{\BL T,T';\, R}_{N}$ up to the boundary---compare with the right part of inequalities \eqref{E:normHs}:

\begin{lemma}
 \label{L:OmR}
  Let $0<  R\le R_{0}$. Then
 for any $s\ge0$ for which $\gu_1\in H^{s}(\Omega_{1})$, there is a constant $C$ such that  
 for all $N\ge1$
\begin{equation}
   \label{E:HsOmR}
 \DNorm{\gu_{N}}{H^{s}(\Omega^{\BL 0,2\pi;\, R}_{N})} 
  \le C\,N^{s-\frac32} \,.
\end{equation}
\end{lemma}

Note that this is typical boundary layer behavior: The ``thickness'' of $\Omega^{\BL 0,2\pi;\, R}_{N}$ is bounded by $1-R^{\frac1N}=O(N^{-1})$. On the whole domain $\Omega_{N}$, the $L^{2}$-norm, for example, cannot behave as $O(N^{-\frac32})$; it is easy to see that on a fixed disk $B_{R'}(0)$, both the $L^{\infty}$ and the $L^{2}$- norm of $\gu_{N}$ are of order $O(N^{-1})$. This is due to the constant term $\gu_{N}(0)=\frac1N\gu_1(0)>0$.

So we need a different argument for the norm on the ball $B_{R^{1/N}}(0)$ which satisfies
\[
   \overline\Omega_N = \overline\Omega^{\BL 0,2\pi;\, R}_{N} \cup \overline B_{R^{1/N}}(0) .
\]
For this, we choose $R_1>0$ such that on $B_{R_1}(0)$, the analytic function $\gu_1$ is given by a convergent Taylor series. 
\begin{equation}
\label{E:Tayloru}
 \gu_1(z) = \sum_{\ell,\ell'\ge0}  a_{\ell\ell'}z^{\ell}\overline z^{\ell'}\qquad (|z|=r< R_1)\,.
\end{equation}
By \eqref{eq:guN},
\begin{equation}
\label{E:TayloruN}
 \gu_{N}(z) = \frac1N\,\sum_{\ell,\ell'\ge0}  a_{\ell\ell'}z^{\ell N}\bar z^{\ell'N}\qquad (|z|=r< R^{\frac1N}_1)\,.
\end{equation}
Integrating this over $B_{R^{1/N}_1}(0)$, one gets 
\begin{equation}
\label{E:d0RN}
  \DNormc{\gu_{N}}{L^2(B_{R^{1/N}_1}(0))}2 \le C\, 
   N^{-2}\,R^{2}_1 \,.
\end{equation}

For a derivative $\partial_{z}^{j_{1}}\partial_{\overline z}^{j_{2}}\gu(z)$,
the convergent series \eqref{E:TayloruN} 
allows immediately to obtain an estimate: For $|j|=m>0$ there is a constant $C_{m}$ independent of $N$
such that
\begin{equation}
\label{E:dmr}
  |\partial^{j}\gu_{N}(z)| \le C_{m}\, N^{m-1}\,r^{N-m} 
  \qquad (|z|=r< R^{\frac1N}_1\,,\; N\ge m) \,.
\end{equation}
Integrating this over $B_{R^{1/N}_1}(0)$, one gets  
with a constant $C$ independent of $N$,
\begin{equation}
\label{E:dmRN}
  \Normc{\gu_{N}}{H^{m}(B_{R^{1/N}_1}(0))}2 \le C\, 
   \tfrac{N^{2m-2}}{2N-2m+2}\,R^{\frac{2N-2m+2}{N}}_1
  \le C_m\, N^{2m-3}\,.
\end{equation}

\begin{lemma}
 \label{L:upper}
 For any integer $m\ge1$ for which $\gu_1\in H^{m}(\Omega_{1})$, there is a constant $C$ such that for all 
 $N\ge m$
\begin{equation}
   \label{E:HmOm}
 \DNorm{\gu_{N}}{H^{m}(\Omega_{N})} 
  \le C\,N^{m-\frac32} \,.
\end{equation}
\end{lemma}
\begin{proof}
Combine Lemma~\ref{L:OmR} for $R<\min(R_0,R_1)$
 with estimates \eqref{E:d0RN} and \eqref{E:dmRN}. 
 \end{proof}

\subsection{Truncation}\label{SS:local}

Recall that we have defined, see \eqref{eq:wtQN}, finite rectangles $\widetilde Q_{N,\rm per}$ and $\widetilde Q_{N}$ by bounding $y$ by $Y_0$ where $Y_0>\max h$. Then $\xi$ is an diffeomorphism from $\widetilde Q_{N,\rm per}$ onto the domain
\[
   \Omega_N\setminus \overline B_{\re^{-Y_0}}(0).
\]

 Recall now the functions $v_N$, $\gu_N$, and $u_N$ from \eqref{E:vN} and \eqref{E:uNdef}, and their counterparts in $\widetilde Q_N$, i.e. $V_N=v_N\circ\xi$, $\gU_N\circ\xi$,  and $U_N=u_N\circ\xi$.
We deduce from Lemma \ref{L:upper} for $m=1$
\begin{equation*}
   \DNorm{\gU_{N}}{H^{1}(\widetilde Q_{N})} 
  \le C\,N^{-\frac12} \,,
\end{equation*}
and we recall from the left part of \eqref{E:estHs} and Lemma \ref{L:nonconst}
\begin{equation*}
   \DNorm{\gU_{N}}{H^{\frac32}(Q^{a,b,\alpha,\beta}_{N})} 
  \ge C>0 ,\qquad \mbox{for}\quad N\ge \tfrac{4\pi}{b-a}\,.
\end{equation*}
Now, for $U_N=u_N\circ\xi=\log|\log\delta_N|\;\gU_N$, we obtain
\begin{equation}
\label{E:wtQNH1}
   \DNorm{U_{N}}{H^{1}(\widetilde Q_{N})} 
  \le C\,N^{-\frac12} \log|\log\delta_N|\,,
\end{equation}
and, for a positive $C$ independent of $N$ 
\begin{equation}
\label{E:wtQNH32}
   \DNorm{U_{N}}{H^{\frac32}(Q^{a,b,\alpha,\beta}_{N})} 
  \ge C\log|\log\delta_N| ,\qquad \mbox{for}\quad N\ge \tfrac{4\pi}{b-a}\,.
\end{equation}
Moreover, from the boundedness of $v_N$ in the $H^{\frac32}(\Omega_{N})$-norm (Lemma \ref{L:veps}) we deduce the boundedness of $V_N$ in the $H^{\frac32}(\widetilde Q_{N})$-norm.
\smallskip

 We choose now numbers $0< a_{1}<a_{2}<b_{2}<b_{1}<2\pi$ and $0<Y_{2}<Y_{1}<Y_{0}$, and define the cut-off function $\chi\in C_{0}^{\infty}(\R\times[0,+\infty))$ as
\[
   \chi(x,y) = \varphi(x)\,\psi(y)
\]
where
\[
\begin{aligned}
   &\mbox{$\varphi\in C_{0}^{\infty}(\R)$ with support in $(a_1,b_1)$ and $\varphi=1$ on $[a_2,b_2]$}\\
   &\mbox{$\psi\in C_{0}^{\infty}([0,+\infty))$ with support in $[0,Y_1)$ and $\psi=1$ on $[0,Y_2]$.}
\end{aligned}
\] 
\smallskip

We define the truncation of $W_N=U_N-V_N$ as
\begin{equation}
\label{eq:wtWN}
  \widetilde W_{N} = \chi W_{N} = \chi (w_N\circ\xi)\,.
\end{equation}
  
 Then the sequence $(\widetilde W_{N})$ is an example for condition (ii) of Theorem~\ref{T:blowup} (hence proving Proposition \ref{P:bc}) in the following sense.

\begin{proposition}
 \label{P:wtilde}
There holds:  
$\widetilde W_{N}\in H^{1}_{0}(\widetilde Q_N)$ and
\begin{itemize}
\item[(a)]\quad $\Norm{\widetilde W_{N}}{H^{\frac32}(\widetilde Q_N)}\to\infty$ \;as \;$N\to\infty$,
\item[(b)]\quad $\DNorm{\Delta\widetilde W_{N}}{H^{-\frac12}(\widetilde Q_N)}$
  \;remains bounded as \;$N\to\infty$.
\end{itemize}
\end{proposition}

\begin{proof}
(a) \;By construction $w_{N}\in H^{1}_{0}(\Omega_{N})$, hence 
$$
\widetilde W_{N}\in H^{1}_{0}(Q_{N}\cap \mathop{\rm supp}\chi)
\subset H^{1}_{0}(\widetilde Q_{N})\,.
$$
Let us choose $\alpha$ and $\beta$ such that $\max h<\alpha<\beta<Y_2$.
On $Q^{a_2,b_2,\alpha,\beta}_{N}$ we have
$\widetilde W_{N}=W_{N}=U_{N}-V_{N}$, hence we can use the estimate \eqref{E:wtQNH32} and the boundedness of $V_N$ in $H^{\frac32}(\widetilde Q_{N})$-norm to get the estimate for $\widetilde W_N$ with constants independent of $N$:
\begin{equation}
\label{E:twlower}
  \Norm{\widetilde W_{N}}{H^{\frac32}(Q^{a_2,b_2,\alpha,\beta}_{N})}
  \ge
  c_{1}\,\log|\log\delta_N| - c_{2}\,.
\end{equation}
Since $Q^{a_2,b_2,\alpha,\beta}_{N}\subset \widetilde Q_{N}$, we see that 
$\Norm{\widetilde W_{N}}{H^{\frac32}(\widetilde Q_{N})}\to\infty$.
\smallskip

\noindent(b) \;For 
$\Delta \widetilde W_{N}= \chi\Delta W_{N} + 2\nabla\chi\cdot\nabla W_{N}
  + W_{N}\Delta\chi$, we have to estimate three terms:
\smallskip

1. For $\chi\Delta W_{N}=-\chi\Delta V_{N}$, we have seen in Section~\ref{SS:rhs} (Lemma \ref{L:fN}) that the $H^{-\frac12}(\D)$-norm of $\Delta v_N$ remains bounded, hence the same for the $H^{-\frac12}(\T\times(0,Y_0))$-norm of $\Delta V_N = \re^{2y}v_N\circ\xi$.
\smallskip

2. For $\nabla\chi\cdot\nabla W_{N}=\nabla\chi\cdot(\nabla U_{N}-\nabla V_{N})$, we use \eqref{E:wtQNH1} which implies that the $H^{1}$-norm of $U_{N}$ still tends to zero, which implies that the $L^{2}$-norm of this term, and a fortiori the $H^{-\frac12}$-norm, remains bounded.
\smallskip

3. For $W_{N}\Delta\chi$, we know that even the $H^{1}$-norms of $U_{N}$ and $V_{N}$ remain bounded.
\end{proof}

\subsection{Conclusions}\label{SS:concQN}
The outcome of the previous subsection shows that a cut-off in slow variables $(x,y)$ does not change properties which arise from the rapidly oscillating structure of $U_N$ related to the rapid variable $Nx$. Proposition \ref{P:wtilde} implies Proposition \ref{P:bc} for the family of domains $(\widetilde Q_N)$, and also for their periodic analogues $\widetilde Q_{N,\rm per}$ (note that for the proof in the periodic rectangles, the cut-off in $y$, $\chi=\psi$, would be enough).
\smallskip

Now, in preparation of the proof of the non-closed range result of Theorem \ref{T:NotH32}, we will glue together scaled versions of different $\widetilde Q_N$. For this, two properties makes the task easier:
\begin{itemize}
\item Homogeneity of semi-norms by scaling,
\item Locality of norms.
\end{itemize}

Apart the case $s=0$, the $H^s$-norms are not homogeneous by scaling, but semi-norms with positive $s$ are.

The locality means that if a function $f$ has its support in a domain $\cO$, then the norm of $f$ in any larger domain $\cO'\supset\cO$ is not enlarged. This property holds for Sobolev norms with positive integer exponents, but not for the fractional order norms. This is why we find advantageous to replace the $H^{-\frac12}$-norms on $\widetilde Q_N$ by the (stronger) weighted $L^2$-norm
\[
   L^2_{\sqrt y}(\widetilde Q_N) = \{f\in L^2_{\sf loc}(\widetilde Q_N)\mid\quad
   (x,y)\mapsto\sqrt y \, f(x,y) \; \mbox{ belongs to }\; L^2(\widetilde Q_N)\}\,.
\]
It is known that if we consider the weight $\sqrt{d_N}$ where $d_N$ is the distance to the boundary of $\widetilde Q_N$, we have $L^2_{\sqrt{d_N}}(\widetilde Q_N)\subset H^{-\frac12}(\widetilde Q_N)$ with bounded embedding. Here we will only use the (easier to prove) result: 

\begin{lemma}
\label{L:L2sqrty}
Let $Q\subset Q'\subset \R\times\R^+$. Let $W\in L^2_{\sqrt y}(Q)$ and set $\widetilde W$ its extension by $0$. Then $\widetilde W$ belongs to $H^{-\frac12}(Q')$ with the estimate for a constant $C$  independent of $Q$, $Q'$, and $W$
\begin{equation}
\label{eq:L2sqrty}
   \DNorm{\widetilde W}{H^{-\frac12}(Q')} \le C \DNorm{W}{L^2_{\sqrt y}(Q)}.
\end{equation}
\end{lemma}

\begin{proof}
The fact that for a function $\phi\in H^{\frac12}(\R^{2})$ with support in $\R\times\R^+$ one has an estimate
$$
  \Normc{\phi}{H^{\frac12}(\R^{2})}2 \ge c\,
  \int_{0}^{\infty}\int_{-\infty}^{\infty} y^{-1} |\phi(x,y)|^{2}dx\,dy 
$$ 
is well known \cite{LionsMagenes68} and easy to prove from the definitions, with 
$c=\int_{R}(1+x^{2})^{-\frac32}dx$.
By duality, this implies \eqref{eq:L2sqrty} for $Q'=\R\times\R^+$ and then for general $Q'$ by monotonicity of the $H^{-\frac12}$ norm with respect to the domain.
\end{proof}

 We have the following quantitative and more precise version of Proposition \ref{P:wtilde}:
\begin{proposition}
 \label{P:wtilde'}
 In the situation of Proposition \ref{P:wtilde} and where we recall that $h$ is the generating profile \eqref{eq:h} of the saws $\widetilde Q_N$, see \eqref{eq:wtQN}, we have the estimates for all $N\ge2$: 
\begin{subequations}
\begin{equation}
\label{eq:wtildea}
   \Norm{\widetilde W_{N}}{H^{\frac32}(\widetilde Q_N)} \ge c_{\rm a} \log\log N
   \quad\mbox{with a constant $c_{\rm a}=c_{\rm a}[h]>0$ independent of $N$}
\end{equation}
and
\begin{equation}
\label{eq:wtildeb}
   \DNorm{\Delta\widetilde W_{N}}{L^{2}_{\sqrt{y}}(\widetilde Q_N)} \le c_{\rm b}
   \quad\mbox{with a universal constant $c_{\rm b}>0$}
\end{equation}
\end{subequations} 
\end{proposition}

\begin{proof}
We revisit the proof of Proposition \ref{P:wtilde}.
\smallskip

\noindent(a) \; We start from the lower bound in \eqref{E:twlower} and have only to prove that 
$c_{1}\,\log |\log\delta_N| - c_{2}$ is estimated from below by $c_{\rm a} \log\log N$. Tracking the estimates leading to \eqref{E:twlower}, we see that $c_2$ does not depend on $h$, while $c_1$ depends on it via Lemma \ref{L:nonconst}. Moreover, according to formula \eqref{eq:deltaN}, we have  $\delta_N=\frac{1}{N}\max h$, leading to 
$|\log\delta_N|\ge \log N - c[h]$, and hence to the desired result. 
\smallskip

\noindent(b) \;
We see that we only have to show that $\Delta V_N$ belongs to $L^{2}_{\sqrt{y}}(\widetilde Q_0)$ with uniformly bounded norms. By construction $V_N(x,y)=A_{\delta_N}(y)$ and $\Delta V_N(x,y)=a'_{\delta_N}(y)$. Since $a'_{\delta_N}(y)$ is $0$ for $y<\delta_N$ and equal to the derivative of $\log|\log y|$ for $\delta_N<y<\frac12$, the result follows from explicit integration:
\[
   y\mapsto (\sqrt y \,\log y)^{-1} \; \in\; L^2((0,\tfrac12)),
\]
independently from the value of $\delta_N$.
\end{proof}

\section{Proof of Theorem \ref{T:NotH32}}\label{S:Conc}

\subsection{A class of Lipschitz domains}

We construct the domain $\cO$ as a subset of the square $(0,1)\times(0,1)$, where the lower boundary $\{(x',0)\mid x'\in(0,1)\}$ is replaced by a rough boundary given by the equation
\[
   x'\mapsto y=H(x')
\]
with rapid oscillations accumulating at the origin.

\smallskip\noindent
The function $H$ depends on three ingredients:
\begin{enumerate}
\item A strictly decreasing sequence $j\mapsto x\sj$ in $(0,1)$ with
\[
   x^{(0)} = 1 \quad\mbox{and}\quad x\sj\to 0 \;\mbox{ as }\; j\to\infty.
\]

\item A non-constant function $h\ge0$ that is Lipschitz $2\pi$-periodic and piecewise $\sC^2$,  such that
\[
   h(0)=h(2\pi)=0 \quad\mbox{and}\quad h\ge1.
\]

\item A sequence of positive integers $N\sj$ such that
\[
   N\sj\to\infty \quad\mbox{as}\quad j\to\infty.
\] 
\end{enumerate}

\noindent
Given these, $H$ is constructed in the following way:

The points $x\sj$ define the intervals $\cI\sj = (x\sj,x^{(j-1)})$ so that $[0,1] = \cup_{j\ge1} \overline\cI\sj$.
The affine transformation
from $(0,2\pi)$ onto $\cI\sj$
is given by
\[
   x \mapsto x\sj + \varepsilon\sj \,x \quad\mbox{with}\quad
   \varepsilon\sj = \frac{x^{(j-1)}-x\sj}{2\pi}\,.
\]
Now $H$ is defined by gluing together the scaled versions $h_{N\sj}:x\mapsto \frac{1}{N\sj} h(N\sj x)$ of the function $h$, cf \eqref{eq:hN}, via the above transformations:
\[
H(x') = 
   \begin{cases}
   \varepsilon\sj h_{N\sj}\big( (\varepsilon\sj)^{-1}(x'-x\sj) \big)
   & \mbox{if }x'\in\cI\sj \\
   0 & \mbox{if }x'\in\{x^{(0)}, x^{(1)},\ldots x\sj, \ldots\}.
\end{cases}
\]
It is clear that $H$ is Lipschitz with the same Lipschitz constant as $h$.
We are going to prove:

\begin{proposition}
\label{P:cO}
Let the domain $\cO$ be constructed as above. Then $\cO$ is an instance for the statement of Theorem \ref{T:NotH32}.
\end{proposition}

\begin{proof}
Denote by $\cH\sj$ the affine transformation in two variables
\[
   \cH\sj: (x,y) \mapsto (x',y')\quad\mbox{with}\quad 
   x' = x\sj + \varepsilon\sj \,x\quad\mbox{and}\quad y' = \varepsilon\sj \,y,
\]
set $\cK\sj$ the inverse of $\cH\sj$, and define the domains $\cO\sj$ as the image by $\cH\sj$ of the bounded subdomains $\widetilde Q_{N\sj}$ (cf \eqref{eq:wtQN} where we choose $Y_0=1/\varepsilon^{(1)}$ for instance)
\[
   \cO\sj = \cH\sj(\widetilde Q_{N\sj}).
\]
The lower boundary $\Gamma\sj$ of $\cO\sj$ is given by
\[
   \cI\sj \ni x'\mapsto \varepsilon\sj h_{N\sj}\big( (\varepsilon\sj)^{-1}(x'-x\sj) \big)
\]
 and the lower boundary of $\cO$ is the union of $\overline\Gamma\sj$ for $j\ge1$. Note that (see Fig. \ref{F:OnotH32})
\[
    \cup_{j\ge1} \overline\cO\sj \subset \overline\cO.
\]

\begin{figure}[h]
\centering
\includegraphics[width=0.48\textwidth]{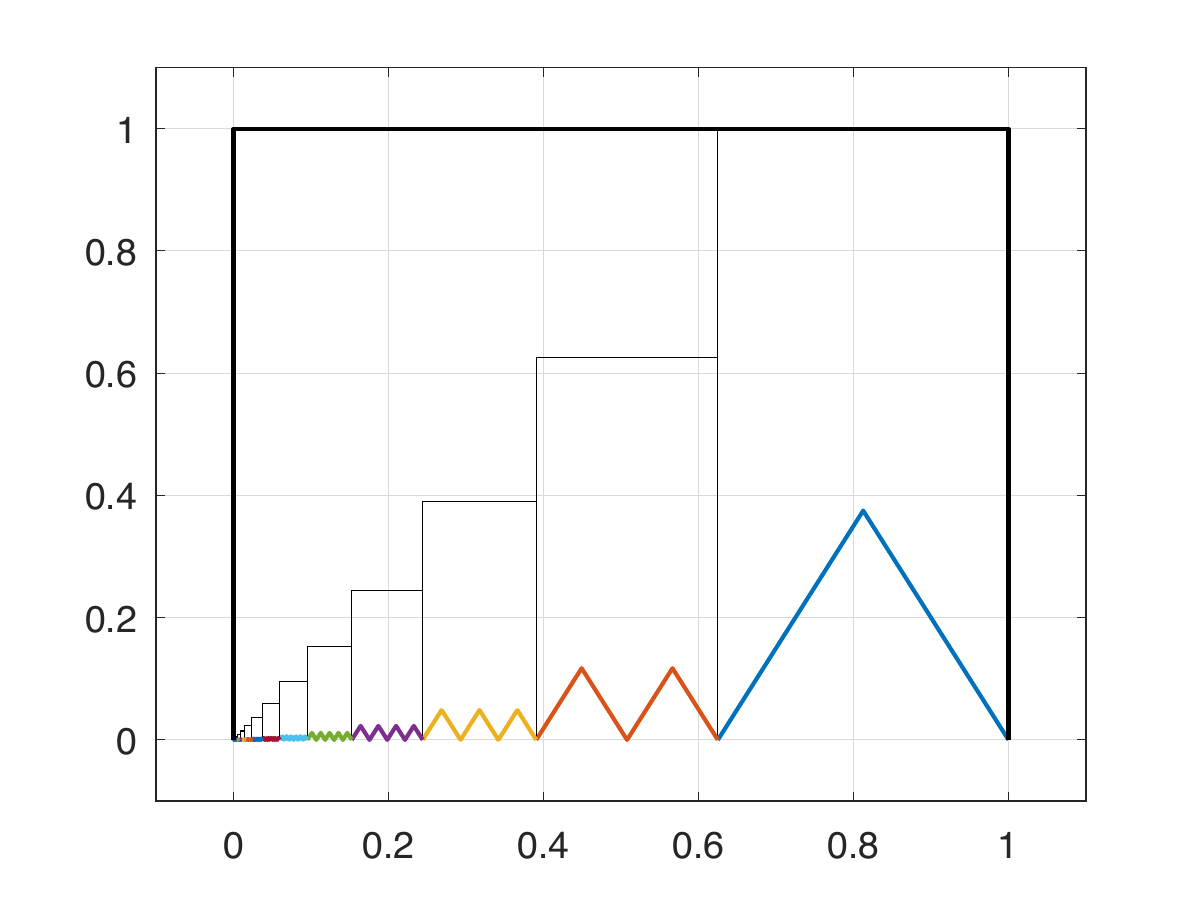}\quad
\includegraphics[width=0.48\textwidth]{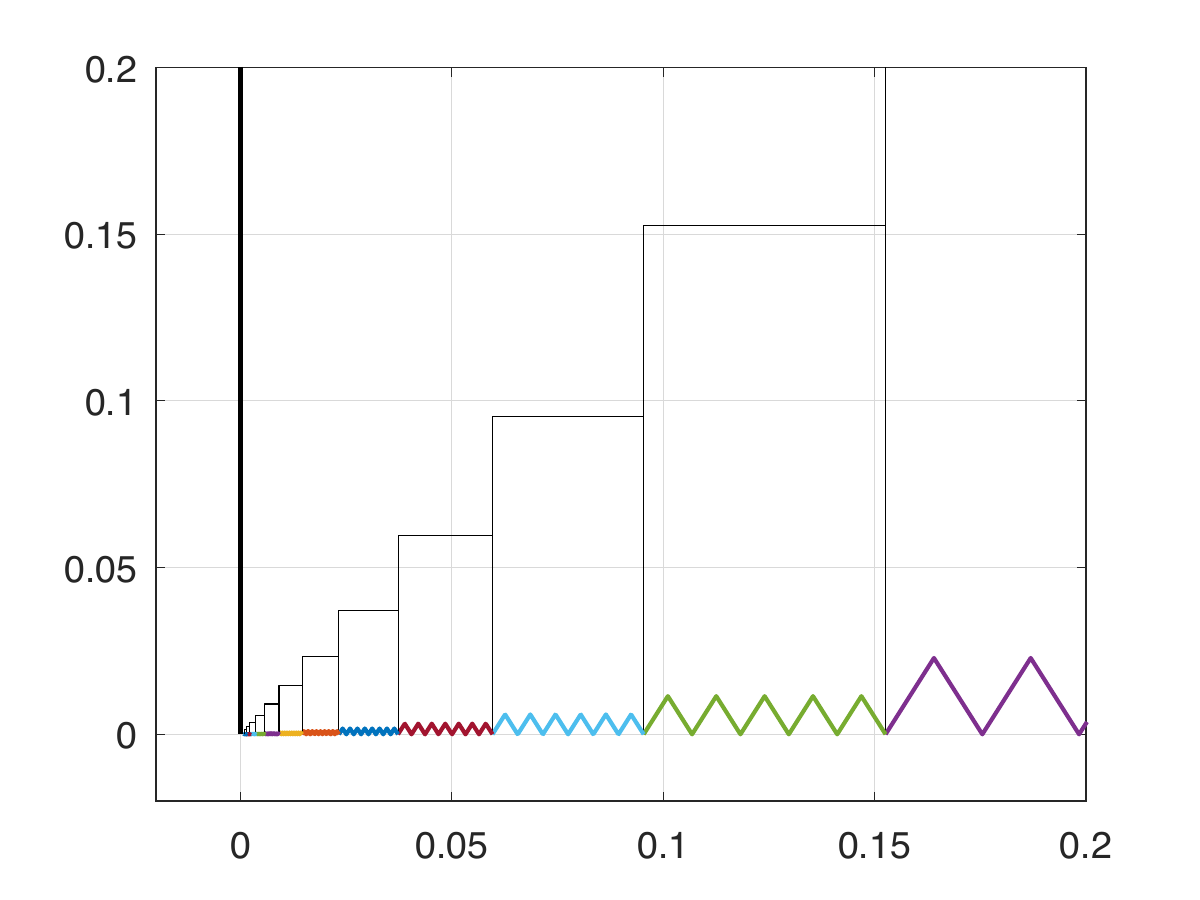}
\caption{\small Domain $\cO$ for $x\sj=1.6^{-j}$, $N\sj=j$, and $h$ the hat function. Zoom on the right. The different colors correspond to the lower boundary above $\cI\sj$. Thin black lines bound the subdomains $\cO\sj$.}
\label{F:OnotH32}
\end{figure}

\noindent
We set, see \eqref{eq:wtWN},
\[
   w\sj = \widetilde W_{N\sj}\circ \cK\sj\quad\mbox{defined on }\quad \cO\sj.
\]
We know from Proposition \ref{P:wtilde'} that the $H^{\frac32}$-seminorms of $\widetilde W_{N\sj}$ blow up as $j\to\infty$ as the sequence $c_{\rm a}[h]\log\log {N\sj}$, while their $H^1$-norm and the $L^2_{\sqrt y}$-norm of $\Delta \widetilde W_{N\sj}$ stay bounded. By the change of variables $\cK\sj$ we find
\[
   \Norm{w\sj}{H^s(\cO\sj)} = (\varepsilon\sj)^{1-s}
   \Norm{\widetilde W_{N\sj}}{H^s(\widetilde Q_{N\sj})},\quad s=0,\;1,\; \tfrac{3}{2}
\] 
and
$\Delta w\sj = (\varepsilon\sj)^{-2} \big(\Delta \widetilde W_{N\sj}\big)\circ\cK\sj$
while $\DNorm{F\circ \cK\sj}{L^2_{\sqrt y'}(\cO\sj)} = (\varepsilon\sj)^{\frac{3}{2}}
 \DNorm{F}{L^2_{\sqrt y}(\widetilde Q_{N\sj})}$. Hence
\[
   \DNorm{\Delta w\sj}{L^2_{\sqrt y'}(\cO\sj)} = 
   (\varepsilon\sj)^{-\frac{1}{2}} 
   \DNorm{\Delta \widetilde W_{N\sj}}{L^2_{\sqrt y}(\widetilde Q_{N\sj})}\,.
\]
Considering now $\breve w\sj$, the extension by $0$ of $w\sj$ from $\cO\sj$ to $\cO$, we see by localization that the norms of integer order are unchanged, whereas the $H^{\frac32}$-seminorm may increase. Therefore, we find for the quotients of norms:
\[
   \frac{\Norm{\breve w\sj}{H^{\frac32}(\cO)}}
   {\DNorm{\Delta \breve w\sj}{L^2_{\sqrt y'}(\cO)}} \ge
   \frac{\Norm{ w\sj}{H^{\frac32}(\cO\sj)}}
   {\DNorm{\Delta  w\sj}{L^2_{\sqrt y'}(\cO\sj)} } 
   =
    \frac{ \Norm{\widetilde W_{N\sj}}{H^{\frac32}(\widetilde Q_{N\sj})}}
   {\DNorm{\Delta \widetilde W_{N\sj}}{L^2_{\sqrt y}(\widetilde Q_{N\sj})} } .
\]
Using Proposition \ref{P:wtilde'} we find
\[
   \frac{ \Norm{\widetilde W_{N\sj}}{H^{\frac32}(\widetilde Q_{N\sj})}}
   {\DNorm{\Delta \widetilde W_{N\sj}}{L^2_{\sqrt y}(\widetilde Q_{N\sj})} } \ge
   \frac{c_{\rm a}[h] \log \log N\sj}{c_{\rm b}}
\]
hence 
\[
   \frac{\Norm{\breve w\sj}{H^{\frac32}(\cO)}}
   {\DNorm{\Delta \breve w\sj}{L^2_{\sqrt y'}(\cO)}}
    \to\infty\quad\mbox{as}\quad j\to\infty.
\]
Moreover, by construction $\breve w\sj$ belongs to $H^1_0(\cO)$. 
Thus there cannot be any a-priori estimate \eqref{E:apriori}, and 
Theorem \ref{T:NotH32} is proved via Proposition \ref{P:cO}.
\end{proof}

\begin{remark}
Taking advantage of the constructions performed in the latter proof it is easy to exhibit a function $w\in H^1_0(\cO)$ such that $\Delta w\in H^{-\frac{1}{2}}(\cO)$ and such that $w\not\in H^{\frac{3}{2}}(\cO)$. It suffices to choose a sequence $(a\sj)_j$ satisfying
\[
   \big(a\sj\big)_j \in \ell^2\quad\mbox{and}\quad \big(a\sj \log\log N\sj\big)_j \not\in\ell^2
\]
and to set
\[
   w = \sum_j a\sj \sqrt{\varepsilon\sj} \;\breve w\sj.
\]
\end{remark}

\subsection{A class of $\sC^1$ domains}
By a modification of the constructions above, we can prove without further effort that Theorem \ref{T:NotH32} still holds if we replace ``Lipschitz" by ``$\sC^1$". For this, instead of the square $\cO_0=(0,1)\times(0,1)$ we start from  a $\sC^2$ domain containing $\cO_0$ and whose boundary contains the segment $[-1,2]\times\{0\}$. We replace the part $[0,1]\times\{0\}$ of this segment by the graph of a function $H$ and revisit the three ingredients producing $H$. 
\begin{enumerate}
\item is unchanged: We take a strictly decreasing sequence $j\mapsto x\sj$ in $(0,1)$ with
\[
   x^{(0)} = 1 \quad\mbox{and}\quad x\sj\to 0 \;\mbox{ as }\; j\to\infty.
\]

\item A non-constant function $h\ge0$ that is $\sC^2$ $2\pi$-periodic,  such that
\[
   h(0)=h(2\pi)=0 \quad\mbox{and}\quad h\ge1,
\]
and (this is the main novelty) we set $h\sj=\frac{1}{j}h$.

\item A sequence of positive integers $N\sj$ satisfying the stronger condition
\[
   c_{\rm a}[h\sj] \log\log N\sj \to\infty \quad\mbox{as}\quad j\to\infty.
\] 
\end{enumerate}
Then $H$ is defined by
\[
   H(x') = \varepsilon\sj (h\sj)_{N\sj}\big( (\varepsilon\sj)^{-1}(x'-x\sj) \big)
   \quad \mbox{for}\quad x'\in\cI\sj.
\]
Here the scaling of the smaller profile $h\sj$ is used instead of $h$ itself. With such a choice, we see that $H'(x')\to0$ as $x'\to0$, proving that the domain $\cO$ is $\sC^1$. Then the functions $w\sj$ are defined on $\cO\sj$ by the same formulas as above:
\[
   w\sj = \widetilde W_{N\sj}\circ \cK\sj
\]
where $\widetilde W_{N\sj}$ has to be understood as stemming from the profile $h\sj$ instead of $h$, whence the condition on $N\sj$ involving $c_{\rm a}[h\sj]$.

\bibliographystyle{siam}
\bibliography{biblio}
\end{document}